\documentclass[12pt,a4paper]{article}%

\usepackage{amsmath,amssymb,amsfonts}

\usepackage{caption2}
\usepackage[]{hyperref}

\usepackage[pdftex]{color,graphicx}
\usepackage{multicol}

\numberwithin{equation}{section} 
\setcaptionmargin{1cm}

\newcommand{\be}{\begin{equation}}
\newcommand{\ee}{\end{equation}}
\newcommand{\m}[1]{{#1}}

\usepackage{color}
\definecolor{rosso}{cmyk}{0,1,1,0.4}
\definecolor{rossos}{cmyk}{0,1,1,0.55}
\definecolor{rossoc}{cmyk}{0,1,1,0.2}

\begin{document}


$\quad$
\vskip 1.0cm
\begin{center}
{\huge \bf \color{rossos} 
An extension to the complex plane of the Riemann-Siegel Z function.
  } 
\vskip 1.0cm {\large

Giovanni  Lodone\footnote{ 
Retired. Email: giolodone3@gmail.com
} 
\\[1cm]
  }

\vskip 1.0cm

\end{center}

\begin{abstract}

 \noindent The usual Riemann-Siegel Z(t) is a real-valued function. We construct a complex function depending from t and from distance from critical line. It is linked to Riemann Xi(s) function by the same real scaling factor of the usual Riemann-Siegel Z(t) on critical line. Errors are not greater than the errors of Riemann-Siegel Z(t) on the critical line, while this result covers at least the whole critical strip.

{

\noindent MSC-Class  :   11M06, 11M26, 11M99

  } 
\vspace{0.1cm}


\end{abstract}


                      \tableofcontents




\section{Introduction}

In 1859 Bernhard Riemann in an outstanding paper  \cite{Riemann:1859ar} \cite[p.~299]{Edwards:1974cz} \cite{Bombieri:2000cz} suggested that the non-trivial
zeros of the analytic continuation of the function:
\be \label {1p1}
\zeta(s)=\sum _{n=1}^\infty1/n^s 
 \quad  \qquad    \qquad \quad \Re(s)>1
\ee
to the so-called critical strip $0<\Re(s)<1 $ must all lie on the line $\Re(s)= 1/2$ (the critical line). The complex argument $s$ is expressed throughout as:
\be \label{1p2}
 s=\frac{ 1}{2}+\epsilon +it
 \ee
so that $\epsilon$ is the distance from critical line parallel to the real axis  and $t$ is the imaginary coordinate.
Until now, this conjecture has remained  unproven, and  is referred to as the Riemann Hypothesis (RH)  \cite{Bombieri:2000cz}. 
A key finding in this effort has been the introduction of the function $\xi(s)$  \cite{Riemann:1859ar}   \cite[p.~16]{Edwards:1974cz}, that is real on the critical line and has the same zeros of $\zeta(s)$ in the critical strip:
\be \label {1p3}
\xi(s)=\Gamma \left( \frac{s}{2}  +1\right)(s-1)\pi^{-s/2}\zeta(s) \, .
\ee
In 1932 C. L. Siegel  crucially succeeded in mastering the information contained  in Riemann's  private notes  about the values of the $\xi$ function on the critical line 
 \cite[p.~136]{Edwards:1974cz} \cite[p.~3]{AleksandarIvic:2003}.    
His main result was the approximate formula:
     \be  \label {SiegelForm} 
     Z(t) =2
     \sum_{n = 1}^{n=N}\frac{  \cos\left(   t  \ \ln\left(    \sqrt{\frac {t}{2  e \pi n^2} } \right) - \frac{\pi}{8}   \right)    }{\sqrt{n}}
+R_0(t) \sim
 - \frac{  \xi( \frac{ 1}{2} +it) } { f(t)}
 \quad ; \quad N=  \left  \lfloor   \sqrt{\frac {t}{2   \pi } },  \right \rfloor 
     \ee
where $\sim$ means   asymptotically equal for  $t    \rightarrow \infty$  and $f(t)$ is the scale factor used in  \cite[p.~176]{Edwards:1974cz}, reported in eq. (\ref  {smallScaleFactor}). 
Equation (\ref {SiegelForm}) has been the core of large scale computations of  the zeros  of $\zeta(s)$ as described  in  \cite {Gourdon:2004cz} and \cite[p.~4]{AleksandarIvic:2003}.

\begin{figure}[!htb]
\begin{center}
\includegraphics[width=0.87\textwidth]{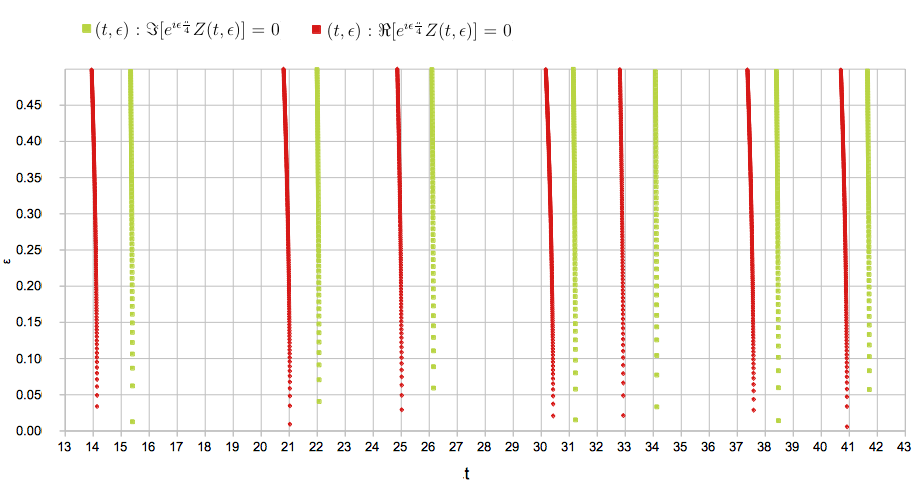} 
\caption{\m{ \small     
The points of  $\Im [ Z(t,\epsilon) ]=0$ and  $\Re [ Z(t,\epsilon) ]=0$  are computed on equispaced  lines at  constant   $t= \Delta t \times m \ : m = 1,2,3 . . .$.  
At the intersection between  the curve  $\Im[ \xi(t,\epsilon) ] \sim \Im [ Z(t,\epsilon) ]=0 $ with the critical line, at $\epsilon =0$, the conformality of the  $\xi(t,\epsilon)$, seen as a complex transformation, is lost. These points are extremal points for the amplitude of  $\xi(1/2+it)$ . At low $t$ values  this is not true  for $Z(t,\epsilon)$  because of the distorting effect of the scale factor $F(t)$.  Note the first  seven zeros at $t= 14.13...;21.02...;25.01...;30.42....;32.93.....;37.58....;40.91....$
} }
\label {T13T42ReImZero}  
\end{center}
\end{figure}

\begin{figure}[!htb]
\begin{center}
\includegraphics[width=0.87\textwidth]{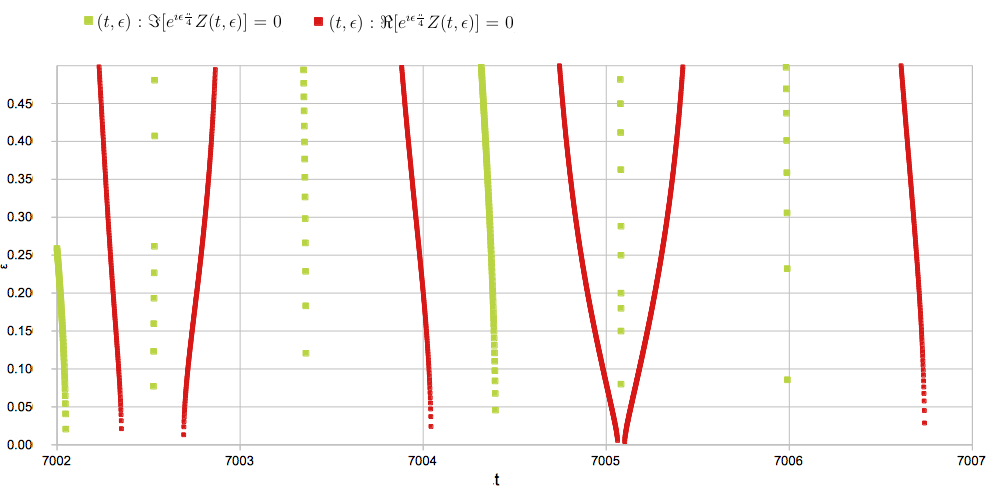} 
\caption{\m{ \small   See       \cite[p.~176]{Edwards:1974cz}          and      \cite[p.~6]{AleksandarIvic:2003}   . 
First Lehemer phenomenon is at  $t \approx$ 7005.  Note the much simpler plot of $\Re[\xi]=0$ and $\Im[\xi]=0$ curves  with respect to the plot of   $\Re[\zeta]=0$ and $\Im[\zeta]=0$  curves  in $(t,\epsilon)$ plane shown   in \cite[p.~342]{ConreyTutorialMarch2003}.
 } }
\label {LehemerOne}
\end{center}
\end{figure}

The purpose of this paper is to develop  an extension  of the Riemann-Siegel function (\ref  {SiegelForm}) that provides a means of expressing  $\xi(\frac{ 1}{2}+\epsilon +it)$ at least   for $ -1<\epsilon  < 1$  using a scale factor defined in eq. (\ref  {FiDiT}):
$$
F(t) =  \left(\pi/2 \right)^{0.25}t^{\frac{7}{4}} e^{- \frac{\pi}{4}  t} 
$$
and  hyperbolic functions.
For ease of reading,  we present our final result, given in eq.  (\ref{ZSinhECosh}):

\be \label {ZSinhECosh1}
 Z(t,\epsilon)=
2 \sum_{n=1}^N \frac{    \cosh \left[ \epsilon \  \ln\left( \sqrt{\frac {t}{2 \pi n^2}  }\right) \right]   }{\sqrt{n}}
\cos \left(   t \  \ln\left(    \sqrt{\frac {t}{2  e \pi n^2} } \right) - \frac{\pi}{8}   \right) + 
\ee

$$  + 2 i  \sum_{n=1}^N \frac{  \sinh\left[ \epsilon \ \ln\left( \sqrt{\frac {t}{2 \pi n^2}  }\right) \right]   }{\sqrt{n}}
\sin \left(   t  \  \ln\left(    \sqrt{\frac {t}{2  e \pi n^2} } \right) - \frac{\pi}{8}  \right )         +R(t,\epsilon)
+ Err(t,\epsilon)  = %
\frac{-\xi( 1/2+\epsilon +it)}{F(t)  e^{i\epsilon\frac{\pi}{4}}}
$$

\noindent where %
$N=  \left  \lfloor   \sqrt{\frac {t}{2   \pi } }  \right \rfloor $,
$F(t) =  \left(\pi/2 \right)^{0.25}t^{\frac{7}{4}} e^{- \frac{\pi}{4}  t}  \ $,  \cite[p.~119]{Edwards:1974cz} for $t$ big, $ |Err(t,\epsilon)| < e^{-0.1 \times t}$ at least for $|t|>100$  ( \cite[p.~144]{Edwards:1974cz},
see  (\ref  {UBL0}),  (\ref  {UBL2}), and,  (\ref  {UBL3}).  $R_1(t)$  is given by eq. (\ref{RQuasiTutto}) specialized  below for $M=1$:

\be
 R_1(t,\epsilon) = (-1)^{N-1} \left( \frac{2 \pi}{t}\right)^{1/4}  \left[     C_0(p)+C_1(p,\epsilon)\left( \frac{2 \pi}{t}\right)^{1/2}   \right]
   \quad ; \quad p= 
   \sqrt{\frac {t}{2 \pi}} -N
\ee
where $ C_0(p)= \cos(2 \pi(p^2-p-1/16) /  \cos(2 \pi p) $.
Using $C_1(p,\epsilon)$ the Z function given in eq.  (\ref{ZSinhECosh1}) is almost  holomorphic for  $t>>1$ and small  $\left| \frac{\epsilon}{t} \right|$, as shown in eq. (\ref  {CauchyRiemannForC1}).         \m{
From a computational point of view it can be meaningful to use eq. (\ref{ZSinhECosh1})  with only $C_0(p)$, that is disregarding $C_1(p,\epsilon)$, since the errors of (\ref{ZSinhECosh1}) with $R_0(t)$ are not higher than the errors of the original Riemann-Siegel (\ref{SiegelForm}), as shown in Figure \ref{MaggiorazioniL023}.

We report in tables \ref{tab1} and \ref{tab2}  the result of a numerical comparison  with {\it Wolfram Mathematica} RiemannSiegelZ function, referred to as $Z_M(t-i\epsilon)$, 
evaluated in same points outside critical line. 

Using  (\ref {ZSinhECosh}) 
 it  is easy, for example, to compute points of a topographical surface defined by  $\Re[e^{i\epsilon\frac{\pi}{4}}  Z(t,\epsilon)]$ above the $(t,\epsilon)$ plane, or to plot in the $(t,\epsilon)$ plane  the  zero-height points  (i.e $(t,\epsilon): \Re[  e^{i\epsilon\frac{\pi}{4}}  Z(t,\epsilon)]=0$).  The same can be done with $\Im[  e^{i\epsilon\frac{\pi}{4}}   Z(t,\epsilon)]$;  see Figures \ref{T13T42ReImZero}  and  \ref{LehemerOne}.
In fact the real or imaginary zero-condition computed for  $e^{i\epsilon\frac{\pi}{4}}   Z(t,\epsilon)$, applies to $\xi(\frac{ 1}{2}+\epsilon +it)$ as well .
The original reason of the present work was precisely to plot such curves for the $\xi$ function, as is done for the $\zeta$ function in  \cite[p.~342]{ConreyTutorialMarch2003}.
 

%
%
}

\begin{table}[bht]
\begin{center}
\begin{tabular}{c|c|c|c|c|c|}
$\epsilon$ & $\Re [Z(t,\epsilon1)]$; see \ref  {ZSinhECosh1} &  $\Im[Z(t,\epsilon)]$ ; see \ref {ZSinhECosh1} & $\Re [Z_M(t- i\epsilon)]$ & $\Im [Z_M(t- i\epsilon)]$ & $||Z-Z_M.||$ \\
\hline
0.1 & 3.241730475804 &- 0.5787043468126&3.241771462370& - 0.578683059836 & 4.62E-005 \\
\hline
0.2 & 3.744381644160 & - 1.199282431530 & 3.744419881769 & - 1.199241919295 & 5.57E-005 \\
\hline
0.3 & 4.642337399179 & -1.907135487575 & 4.642369954238 & -1.907080635877 & 6.38E-005\\
\hline
0.4 & 6.033267420500 & -2.755039800815& 6.033289474138& -2.754980028392 & 6.37E-005 \\
\hline
0.5 & 8.069550741786 & -3.807665184897 & 8.069554328448 & -3.807617504896 & 4.78E-005 \\
\hline
\end{tabular}
\end{center}
\caption{\small Comparison at $ t=7000$  between $Z(t,\epsilon )$ (with only  $C_0(p)$)  in  \ref  {ZSinhECosh1} and {\it Wolfram Mathematica} RiemannSiegelZ function $Z_M(t-i\epsilon)$. Both are evaluated  with a precision of 16 digits. }
\label {tab1}
\end{table}

\begin{table}
\begin{center}
\begin{tabular}{c|c|c|c|c|c|}
$\epsilon$ & $\Re [Z(t,\epsilon)]$; see \ref  {ZSinhECosh1} &  $\Im [Z(t,\epsilon)]$; see \ref  {ZSinhECosh1} & $\Re[ Z_M(t-i\epsilon)]$ & $\Im [Z_M(t- i\epsilon)]$ &  $||Z-Z_M.||$ \\
\hline
0.1 &-0.9050244263086 &0.1402585183494 &-0.9050238453328&0.1402587799723 & 6.37E-007 \\
\hline
0.2 & -1.3080455298552 & 0.2878855199489 & -1.3080449625265 &0.2878860590214 & 7.83E-007\\
\hline
0.3 & -2.143656564420 & 0.441671103478 & -2.143656009082 & 0.441671983111 & 1.04E-006 \\
\hline
0.4 & -3.715336354275 & 0.578566730780 & -3.715335798319 &0.578568148509 & 1.52E-006 \\
\hline
0.5 &-6.582066395361 & 0.629253718227& -6.582065829497 & 0.629256211048 & 2.56E-006 \\
\hline
\end{tabular}
\end{center}
\caption{\small Comparison at $ t=250000$ between $Z(t,\epsilon )$ (with only  $C_0(p)$) in  \ref{ZSinhECosh1} and {\it Wolfram Mathematica} RiemannSiegelZ function $Z_M(t- i\epsilon)$. Both  are evaluated  with a precision of 16 digits.}
\label {tab2}
\end{table}

\noindent The paper is an exercice on  \cite[p.~136-155]{Edwards:1974cz}. The reading will be easier with a copy of chapter 7 of  \cite  {Edwards:1974cz} within reach . 



\section{
Starting point}

In  \cite{Riemann:1859ar}  \cite[p.~137]{Edwards:1974cz}  a contour integral  is given for $\zeta(s) ,\, s\ne 1 , \, s  \in \mathbb{C}$:
\be \label {ZetaContourEdw137}
\zeta(s)  = \frac {\Gamma(1-s ) }{ 2 \pi i} \int_{C} \frac { (-x)^s dx}{(e^x-1) x}
\ee
where $C$ is the contour described in the positive sense starting at $ +\infty$, encircling the origin and returning to $ +\infty $  without crossing the positive real axis.
Afterwards this path  $C$ is deformed continuously in    $C_{N}$, encircling  2N+1  poles of the integral in eq. (\ref {ZetaContourEdw137}).
Evaluating the integral and putting the result in  the definition of $\xi(s)$ in eq. (\ref{1p3}), one finds (see \cite[p.~138]{Edwards:1974cz}, \cite[p.~20]{PUGH:1992ar}):
\begin{eqnarray}
 -\xi(s) &=& (1-s)\Gamma \left(  \frac{s}{2} +1\right)\pi^{\frac{-s}{2}} \left( \sum_{n=1}^{n=N}n ^{-s}  \right) +
(+s)\Gamma\left(  \frac{1-s}{2} +1\right)\pi^{-\frac{1-s}{2}} \left( \sum_{n=1}^{n=N}n ^{-(1-s)}  \right)  +  \nonumber  \\
&+& \frac{ (+s)\Gamma\left(  \frac{1-s}{2} +1 \right)\pi^{\frac{-(1-s)}{2}}   }{  (2\pi)^{s-1} \,  2 \sin(\pi s/2) \, 2\pi i   }\int_{C_{N}} \frac{(-x)^{s-1}e^{-Nx}dx}{e^x-1} \label {Edward(5)pag138} 
\end{eqnarray}
where the sign has been changed for later convenience 
and:
\be\label{app8 2}
N=  \left  \lfloor   \sqrt{\frac {t}{2   \pi } }  \right \rfloor \, .
\ee 

\noindent We now denote by $L_N$ is the usual broken line  whose  path segments $L_0,L_1, L_2$ and $L_3$ are defined in \cite[p.~138]{Edwards:1974cz} \cite[p.~20]{PUGH:1992ar} in order to apply the steepest descent method.  
In Appendix \ref{appendix2} it  is shown that the dominant contribution  to  the integral in          (\ref {Edward(5)pag138}) comes from $L_1$  which  extends from $ a+\frac{1}{2} e^{   i\frac{\pi}{4} |a|  }$ to $ a-\frac{1}{2} e^{   i\frac{\pi}{4} |a|  }$, see (\ref {BoundaryL0L1})  and  (\ref {BoundaryL1L2}), where:
 \be \label {a}
a= i \sqrt{2 \pi t} 
\ee
 is the  saddle point for the evaluation of main integral in (\ref{Edward(5)pag138}), see  Appendix \ref{appendix1} for the details.
We also use the notation of \cite[p.~139]{Edwards:1974cz}, where:
 $$
e^{\phi(x)} =( -x)^{s-1}e^{-Nx}
$$
%
The saddle point occurs when $\phi'(x)=0$ with $\phi(x)= \Re[ (s-1) \  \ln (-x) -Nx ]$ to yield:
\be\label {saddlepointapprox}
\alpha =\frac{  -\frac{1}{2} +it } {N} \approx  \frac{  -\frac{1}{2} +it } {    \sqrt{\frac{t}{2\pi}}   }    \approx 2 \pi i  N  \approx  i \sqrt       {2 \pi t} =a  \quad \mbox{ for } \quad t >> 1
\ee
so that $a$ in (\ref {a}), is an approximation  for the value  $\alpha$  used in \cite[p.~140]{Edwards:1974cz}.
Notice that in (\ref {saddlepointapprox}), $N$ comes from the numerator in $\frac{  -\frac{1}{2} +it } {N} $,  while  in $2 \pi i  N  $ it comes from intergrand denominator ( i.e.  $e^x-1$).

Nothing  changes  in the saddle point integral evaluation procedure with the replacement:
$$
-\frac{1}{2} +it  \quad \rightarrow \quad -\frac{1}{2}+\epsilon +it \quad \mbox{ with } 0 < \epsilon < 1 
$$
as shown in Appendix \ref{appendix1}.
Notice that for $\epsilon = \frac{1}{2} $   the approximation  (\ref{saddlepointapprox}) becomes  exact and we have $\alpha = a $. 


We now make use of the Stirling series  \cite[p.~30]{PUGH:1999ar} that allows us to write:
\be\label {Pi&Gamma2}
\ln( \Gamma(z+1) )=\ln\left( e^{-z}z^{   z+\frac{1}{2}   }(2\pi)^{\frac{1}{2}}  \right )  +\left( \sum_{k=1}^{K-1} \frac{B_{2k}   }{  2 k (2k-1) z^{2k-1}} \right) +R_{2K}(z)
\ee
Although the expression (\ref {Pi&Gamma2})  is a non-convergent asymptotic expansion, it  can be used to estimate the size of the error  $|R_K(z)|$ \cite[p.~112]{Edwards:1974cz}.                 
The $B_i$  are the Bernoulli numbers that vanish for odd $i$ while (see for example  \cite[p.~114]{MurrayRSpiegel:2004cz}):
$$
B_{2} =\frac{1}{6}
\quad ; \quad
 B_{4} =-\frac{1}{30}
\quad ; \quad
B_{6} =\frac{1}{42}
\quad ; \quad
B_{8} =-\frac{1}{30}
\quad ; \quad
B_{10} =\frac{5}{66} 
\quad ; \quad
 B_{12} =-\frac{691}{2730}
\quad ...
$$




        
 
The modulus of  the error term $|R_{2K}(z)|$ is bounded by: 
\be \label  {RestoInApp}
 |R_{2K}(z)| < \left(    \frac{B_{2K}   }{  2 K (2k-1) z^{2K-1}}          \right)
\frac{1}{       \left[     \cos\left ( \frac{\arg(z)}{2}\right)         \right]^{2K}        }
\ee
where $\arg(z)$ is taken in the interval: $-\pi <\arg(z) <\pi$ (see  \cite[p.~112]{Edwards:1974cz} and theorem 2.3 in  \cite[p.~40]{PUGH:1999ar}  due to Stieltjes).
 
 

   
     
       
    
       
       
 


\section{Splitting in simple pieces} \label {Splitting in simple pieces}

To evaluate the expression in (\ref{Edward(5)pag138}) it is convenient to rewrite $s$ in terms of the variable $z$defined by:
     \be \label{sost1}
     z=\frac{s}{2}= \frac{1+2 \epsilon}{4} +\frac{i t}{2} 
     \ee
in the first sum,  so that we have:
$$
     \pi^{\frac{-s}{2}}(1-s) =\pi^{- z}(1-2z) \, ,
$$
while in the second sum and in the third term we use:
\be \label{sost2}
     z= \frac{1-s}{2}= \frac{1-2 \epsilon}{4} -\frac{i t}{2} 
     \ee
 so that:
$$
\pi^{\frac{-(1-s)  }{2}}(+ s) =\pi^{-z}(1-2z)
$$

\noindent In order to manage the computation and also to exploit intermediate results,  it has been devised the notation $A_p^s$  where:

 '' A `` can be:  $A=\Re $ or  $A= \Im$), 
    
``s= +''  means substitution (\ref{sost1})  while  ``s= -''  means substitution (\ref{sost2}).
    
``p'' is for  part (p=1,2 or 3)   
      
Part 1  ( i.e  $\Re^\pm_1$ or  $\Im^\pm_1$ )   is the asymptotic part (as $t \rightarrow \infty$)  of  $\ln\left( e^{-z}z^{   z+\frac{1}{2}   }(2\pi)^{\frac{1}{2}}  \right )$ 
      
Part 2   ( i.e  $\Re^\pm_2$ or  $\Im^\pm_2$ ) refers always to  $\ln\left( e^{-z}z^{   z+\frac{1}{2}   }(2\pi)^{\frac{1}{2}}  \right )$, but as $t \rightarrow \infty$ goes to zero
      
Part 3  ( i.e  $\Re^\pm_3$ or  $\Im^\pm_3$ )  refers to:
      \be \label {BerSumMinusErr} 
      \left( \sum_{k=1}^{K-1} \frac{B_{2k}   }{  2 k (2k-1) z^{2k-1}} \right) 
      \ee

For example for $K=3$, 
setting $\arg(z) \approx \frac{\pi}{2}$ and $B_{2K} = 1/42 $  \m{as  first erased term in (\ref {Pi&Gamma2}),}  the error  $|R_{2K}(z)|$  is bounded by
 (\ref{RestoInApp}) and  \m{using (\ref {sost1}) and (\ref {sost2})  we have}, for the error term:
\begin{eqnarray} 
 |R_{2K}(z)^\pm| &=& |Re_{error}^{\pm} +i \quad Im_{error}^{\pm}| =  \left|R_6\left(\frac{1/2 \pm \epsilon  \pm i t}{2}\right) \right|
\label {MaggErrSerAsint1} \\
&<& \left|  \frac{B_6}{6 * 5 * \left(\frac{1/2\pm \epsilon  \pm i t}{2}\right)^5  *\cos\left (\arg\left(\frac{1/2\pm \epsilon \pm i t}{2}  \right) \right)^6}    \right|  < \left|  \frac{1}{157 \left(\frac{1/2 \pm \epsilon  \pm i t}{2}\right)^5 }  \right|
\nonumber \\
& \le &  \left|  \frac{1}{4.9 t ^5 }  \right|      
\nonumber
\end{eqnarray}
%
It is possible to improve the precision by increasing $K$  but, \m{as we will see in the following
}, if we content ourselves with an error  $< t^{-1}$ 
then  we can completely ignore  the  Bernoulli sum in (\ref{Pi&Gamma2}).

\subsection { First substitution}

If we take the first case (\ref{sost1}), to be used in first sum of (\ref{Edward(5)pag138}), we have  for (\ref  {Pi&Gamma2}):

$$
\ln\left( e^{-z}z^{   z+\frac{1}{2}   }(2\pi)^{\frac{1}{2}} \pi^{- z}(1-2z) \right )_{z= \frac{1+2 \epsilon}{4} +\frac{i t}{2}} =\Re_1^+ +i \Im_1^+ + \Re_2^+ +i \Im_2^+  
$$

$$
 =
\left( \frac{3+2 \epsilon}{4} +\frac{it}{2} \right) \left[  \ln \sqrt{\left( \frac{1+2 \epsilon}{4} \right)^2+\frac{t^2}{4}} + i \quad \arctan\left (         \frac{2t}{1+2\epsilon} \right)  \right]   -
$$

\be \label{gammaS1}
- \left(   \frac{1+2 \epsilon}{4} +\frac{it}{2}     \right) +\frac{\ln(2 \pi)  }{2}  -\frac{\frac{1}{2}+\epsilon+it}{2}  \ln(\pi)+\ln \sqrt{\left( \epsilon-\frac{1}{2}  \right)^2+t^2}  +  i \quad \arctan\left (         \frac{-t}{+1/2-\epsilon} \right) 
\ee

\noindent           We have :

\be \label{TangManip}  
  \arctan\left (         \frac{-t}{+1/2-\epsilon} \right)  = -\frac{\pi}{2} +\arctan\left (         \frac{1-2\epsilon}{2t} \right) \quad and  \quad
  \arctan\left (         \frac{2t}{1+2\epsilon} \right)  = \frac{\pi}{2} -\arctan\left (         \frac{1+2\epsilon}{2t} \right) 
\ee

\noindent             so separating real part from imaginary part we have:

$$
\Re\left( \ln\left( e^{-z}z^{   z+\frac{1}{2}   }(2\pi)^{\frac{1}{2}}\pi^{- z}(1-2z)  \right )_{z= \frac{1+2 \epsilon}{4} +\frac{i t}{2}}  \right) = \Re_1^+ +\Re_2^+ =
$$

$$
=\left( \frac{3+2 \epsilon}{4}  \right)    \left[ \ln\left(\frac{t}{2} \right)+ \frac{1}{2} \ln \left( 1 +\left( \frac{1+2 \epsilon}{2t} \right)^2 \right)   \right]-
$$
$$
- \frac{t}{2} \left[  \frac{\pi}{2} -\arctan   \left (         \frac{1+2\epsilon}{2t} \right)       \right] - \left( \frac{1+2 \epsilon}{4}  \right)+\frac{\ln(2\pi)}{2}-\frac{\ln(\pi)(1+2\epsilon)}{4} + \ln(t) +\frac{1}{2}\ln\left[    1+ \left(   \frac{1-2 \epsilon}{2t}    \right)^2 \right]
$$

\noindent         For  $t>>\epsilon$  we have 
\be \label {Trick1+}
- \frac{t}{2} \left[  -\arctan   \left (         \frac{1+2\epsilon}{2t} \right)       \right] = \frac{1+2\epsilon}{4}  \quad \frac {2t}{1+2\epsilon}\arctan   \left (         \frac{1+2\epsilon}{2t} \right)  
 \quad \rightarrow \frac{1}{4} + \frac {\epsilon}{2}  ;\quad t \rightarrow \infty
\ee

\noindent         And we can write  the asymptotic part not leading to zero as:

\be \label {Trick2+}
- \frac{t}{2} \left[  -\arctan   \left (         \frac{1+2\epsilon}{2t} \right)       \right] =
 \frac{1+2\epsilon}{4} \left[  \quad \frac {2t}{1+2\epsilon}\arctan   \left (         \frac{1+2\epsilon}{2t} \right)   +1 -1\right] = \frac{1+2\epsilon}{4}  
\ee

 \noindent       And we  put it in first piece  $\Re_1^+$ .
 \noindent While:

$$
 \frac{1+2\epsilon}{4} \left[  \quad \frac {2t}{1+2\epsilon}\arctan   \left (         \frac{1+2\epsilon}{2t} \right)    -1\right]
$$
is put in second piece (i.e. $\Re_2^+$), because as  $t\rightarrow \infty $ it goes to zero .

\noindent      Using the notation  $A_p^s$ 
we try  afterward to put together all pieces.

\noindent        $\Re_2^+ + \Re_3^+$ and  $\Im_2^- + \Im_3^-$  
  are developped in powers of $t^{-1}   $,$t^{-2}   $,$t^{-3}   $ . . . in  
  subsection \ref  {Error  merging   }; see (\ref   {Im2e3Pot1}).

\subsection { Computation of $ \Re_1^+  $}

So  $\Re_1^+$ first piece of Real part  with first substitution (\ref{sost1}) is:
$$
 \Re_1^+=  \frac{3}{4}\ln\left(  \frac{t}{2} \right) -\frac{1}{4}+\frac{\ln(2\pi)}{2}-\frac{\ln(\pi)}{4} + \frac{1}{4} +\frac {\epsilon}{2}++\frac{\epsilon}{2} \left(    \ln  \frac{t}{2}   -\ln(\pi)   \right)-
$$

\be\label{asintoticaper gamma1re}
-\frac{\epsilon}{2}-\frac{\pi}{4}t+\ln(t) = \left(
\frac{3}{4}Ln\left(  \frac{t}{2 } \right) +\frac{\ln(2\pi)}{2}-\frac{\ln(\pi)}{4}  -\frac{\pi}{4}t+\ln(t)
\right)
+   \ln \left( \sqrt{ \frac{t}{2 \pi} }   \right) ^\epsilon 
\ee

\subsection { Computation of $ \Re_2^+  $}

While  $\Re_2^+$, second piece of Real part for  first  substitution (\ref{sost1}), is:

$$
\Re_2^+= 
\left( \frac{3+2 \epsilon}{4}  \right)    \left[  \frac{1}{2} \ln \left( 1 +\left( \frac{1+2 \epsilon}{2t} \right)^2 \right)   \right]+
$$

\be\label{SecondoPezzoReGa1}
+\frac{1+2\epsilon}{4} \left[  \quad \frac {2t}{1+2\epsilon}\arctan   \left (         \frac{1+2\epsilon}{2t} \right)    -1\right]
 +\frac{1}{2}\ln\left[    1+ \left(   \frac{1-2 \epsilon}{2t}    \right)^2 \right]   =
\ee

$$
=\left( \frac{3+2 \epsilon}{4}  \right) 
 \left[   
    \frac{1}{2} \left( \frac{1+2 \epsilon}{2t} \right)^2  -
   \frac{1}{4} \left( \frac{1+2 \epsilon}{2t} \right)^4   + . . .  
     \right]  +
     \frac{1+2\epsilon}{4} \left[ 1 -    \frac{1}{3}   \left (         \frac{1+2\epsilon}{2t} \right) ^2 -1+ . . .      \right]+
$$

$$
+\frac{1}{2}\left[     \left(   \frac{1-2 \epsilon}{2t}    \right)^2   -\frac{1}{2}   \left(   \frac{1-2 \epsilon}{2t}    \right)^4 + . . .  \right] =
$$

\noindent  Note that   the most significant powers are    $t^{-2}$ so that, up to the terms of order $t^{-2} $

$$
    \Re_2^+ =\left( \frac{3+2 \epsilon}{4}  \right) 
 \left[   
    \frac{1}{2} \left( \frac{1+2 \epsilon}{2t} \right)^2  - . . .  
     \right]  +
     \frac{1+2\epsilon}{4} \left[ -    \frac{1}{3}   \left (         \frac{1+2\epsilon}{2t} \right) ^2 + . . .      \right]+\frac{1}{2}\left[     \left(   \frac{1-2 \epsilon}{2t}    \right)^2   + . . .   \right]  
     $$

 \noindent It is apparent that all three terms of $\Re_2^+$ go to zero as $t \rightarrow \infty$.

\subsection { Computation of $ \Im_1^+  $}

 For the imaginary part we have:

$$
\Im\left(   ln\left( e^{-z}z^{   z+\frac{1}{2}   }(2\pi)^{\frac{1}{2}} \pi^{- z}(1-2z) \right )_{z= \frac{1+2 \epsilon}{4} +\frac{i t}{2}} \right) = \Im_1^+ + \Im_2^+ =
$$
 
$$
=
\left( \frac{t}{2}  \right)    \left[ \ln\left(\frac{t}{2} \right)+ \frac{1}{2} \ln \left( 1 +\left( \frac{1+2 \epsilon}{2t} \right)^2 \right)   \right] +
$$

\be \label{ParteImag1}
+\left( \frac{3+2 \epsilon}{4}  \right) \left[  \frac{\pi}{2} -\arctan   \left (         \frac{1+2\epsilon}{2t} \right)       \right] -\frac{t}{2}  -\frac{t}{2} \ln(\pi)  -\frac{\pi}{2} +\arctan\left (         \frac{1-2\epsilon}{2t} \right)
\ee

\noindent      For  $t>>\epsilon$

\be\label{asintoticaper gamma1im}
 \Im_1^+=
    \frac{t}{2}\ln\frac{t}{2\pi}  -\frac{t}{2} -\frac{\pi}{8}+\frac{\pi}{4}\epsilon \quad ; \quad t >>\epsilon
\ee

\subsection { Computation of $\Im_2^+  $}

\be \label{ IM2Sost1}
\Im_2^+=     \frac{t}{4} \ln \left( 1 +\left( \frac{1+2 \epsilon}{2t} \right)^2 \right)   -\left( \frac{3+2 \epsilon}{4}  \right) \left[  \arctan   \left (         \frac{1+2\epsilon}{2t} \right)       \right]  +\arctan\left (         \frac{1-2\epsilon}{2t} \right) =
\ee

$$
= \frac{t}{4}  \left( \left( \frac{1+2 \epsilon}{2t} + . . . . \right)^2 \right)-\left( \frac{3+2 \epsilon}{4}  \right) \left[     \left (         \frac{1+2\epsilon}{2t} \right)  +. . . .     \right]  +\left (         \frac{1-2\epsilon}{2t} \right) 
$$
up to terms of order $t^{-1}$. Notice that   $\Im_2^+$ for   goes to zero for $t \rightarrow \infty$.

\subsection { Computation of $ \Re_3^+  $ and  $ \Im_3^+   $}

Let us take now the sum in (\ref{Pi&Gamma2}), always with substitution (\ref{sost1}), for K=3.

\noindent      Here too it is useful to consider real part:

$$\Re_3^+  (only \quad till \quad t^{-2} \quad terms)
  = Re\left( \sum_{k=1}^{K-1} \frac{B_{2k}   }{  2 k (2k-1) z^{2k-1}} \right)_{z= \frac{1+2 \epsilon}{4} +\frac{i t}{2}   \quad ; \quad K=3 } =
$$

\be \label{tillB4RE}
=\frac{1+2 \epsilon}{48\left [        1+ \left( \frac{1+2 \epsilon}{2 t}  \right)^2 \right]    \left( \frac{t}{2} \right) ^2} 
- \frac{  (1+2 \epsilon)^3} {360 t^6  \left(  1+ \left( \frac{1+2 \epsilon}{2 t}  \right)^2   \right)^{3}  }
+\frac{ (1+2 \epsilon)}{30 t^4   \left(  1+ \left( \frac{1+2 \epsilon}{2 t}  \right)^2   \right)^{3}   } =
\ee

$$
=\frac{1+2\epsilon}{12 t^2}-...
$$
%
%
%
and imaginary part:
$$
\Im_3^+  (only \quad till\quad t^{-1} \quad terms) =
 \Im\left( \sum_{k=1}^{K-1} \frac{B_{2k}   }{  2 k (2k-1) z^{2k-1}} \right)_{z= \frac{1+2 \epsilon}{4} +\frac{i t}{2}   \quad ; \quad K=3 } = 
$$

\be \label{tillB4IM}
= - \frac{1}{6 t  \left(  1+\left(   \frac{1+2 \epsilon}{2 t}  \right)^2   \right) }
- \frac{1}{45 t ^3 \left(  1+\left(   \frac{1+2 \epsilon}{2 t}  \right)^2   \right)^3 }
+\frac{(1+2 \epsilon)^2}{60 t ^5 \left(  1+\left(   \frac{1+2 \epsilon}{2 t}  \right)^2   \right)^3} =
\ee

$$
=- \frac{1}{6t}+. . .  
$$

\subsection { Second substitution}
Second substitution  (\ref{sost2})  applies to second sum ant remainder (i.e. 3rd term) in (\ref{Edward(5)pag138}) . 

\noindent      So we have, likewise:

$$
\ln\left( e^{-z}z^{   z+\frac{1}{2}   }(2\pi)^{\frac{1}{2}} \pi^{- z}(1-2z) \right )_{z= \frac{1-2 \epsilon}{4} -\frac{i t}{2}}= \Re_1^- +\Re_2^- +i \Im_1^- + i \Im_2^- =
$$

$$
 =\left(    \frac{3}{4} -\frac{\epsilon}{2}-\frac{it}{2}\right)     \left[ \ln \sqrt{\left (  \frac{1}{4}-\frac{\epsilon}{2} \right)^2+\frac{t^2}{4}    }  + i \,\, \arctan\left (         \frac{-t}{1/2-\epsilon} \right) \right]  +
$$

$$
- \left(   \frac{1-2 \epsilon}{4} - \frac{it}{2}     \right) +\frac{\ln(2 \pi)  }{2} -
$$

$$
 -\frac{\frac{1}{2}-\epsilon-it}{2}\ln(\pi) +\ln \sqrt{\left( \epsilon+\frac{1}{2}  \right)^2+t^2}   \  \ +  i \ \arctan\left (         \frac{+t}{+1/2+\epsilon} \right) 
$$

\noindent    Referring to (\ref{TangManip}),  we separate again the real from the imaginary part.

We find:

$$
\Re\left(\ln\left( e^{-z}z^{   z+\frac{1}{2}   }(2\pi)^{\frac{1}{2}} \pi^{- z}(1-2z) \right )_{z= \frac{1-2 \epsilon}{4} -\frac{i t}{2}} \right)= \Re_1^-+\Re_2^-=
$$

$$
=\left ( \frac{3}{4}  - \frac{\epsilon}{2}\right) \left[    \ln\left(  \frac{t}{2}\right)  +\frac{1}{2}  \ln\left(   1+  \left[    \frac{1+2\epsilon}{2t} \right]^2  \right)  \right] +
$$
$$
+\frac{t}{2} \left[  -\frac{\pi}{2} +\arctan\left (         \frac{1-2\epsilon}{2t} \right)\right]-\left ( \frac{1}{4}  - \frac{\epsilon}{2}\right) +\frac{\ln(2 \pi)}{2} - \left(\frac{1}{2} -\epsilon \right) \frac{\ln(\pi)}{2} +
$$
$$
+ \ln(t) + \frac{1}{2} \ln\left (  1+\left[  \frac{1+2\epsilon}{2t} \right] ^2 \right)
$$

\noindent       For  $t>>\epsilon$  we have 

\be \label {Trick1-}
 \frac{t}{2} \left[  \arctan   \left (         \frac{1- 2\epsilon}{2t} \right)       \right] = \frac{1-2\epsilon}{4}  \quad \frac {2t}{1-2\epsilon}\arctan   \left (         \frac{1- 2\epsilon}{2t} \right)  
 \quad \rightarrow \frac{1}{4}  -  \frac {\epsilon}{2}  ;\quad t \rightarrow \infty
\ee

\noindent         And we can write  the asymptotic part not leading to zero (contained in   $\Re_1^-$) as:


\be \label {Trick2-}
 \frac{t}{2} \left[  \arctan   \left (         \frac{1-2\epsilon}{2t} \right)       \right] =
 \frac{1-2\epsilon}{4} \left[  \quad \frac {2t}{1-2\epsilon}\arctan   \left (         \frac{1-2\epsilon}{2t} \right)   +1 -1\right] = \frac{1-2\epsilon}{4}  
\ee

\noindent       and

$$
\frac{1-2\epsilon}{4} \left[  \quad \frac {2t}{1-2\epsilon}\arctan   \left (         \frac{1-2\epsilon}{2t} \right)    -1\right]
$$
\noindent        which is contained in  $\Re_2^-$,  because as  $t\rightarrow \infty $ it goes to zero.
Thus for  $t >>\epsilon$  and  $t>>1$:
 
 $$
 \Re\left( \ln\left( e^{-z}z^{   z+\frac{1}{2}   }(2\pi)^{\frac{1}{2}} \pi^{- z}(1-2z) \right )_{z= \frac{1-2 \epsilon}{4} -\frac{i t}{2}}  \right)=\Re_1^- + \Re_2^-
$$

\subsection { Computation of $  \Re_1^-   $}

\noindent       We can write:

$$
 \Re_1^-= \frac{3}{4}\ln\left(  \frac{t}{2} \right)-\frac{1}{4}+\frac{1}{4}+\frac{\ln(2\pi)}{2}-\frac{\ln(\pi)}{4}
+\frac{\epsilon}{2} \left(  \ln(\pi) - \ln  \frac{t}{2}    \right)+\frac{\epsilon}{2}-\frac{\pi}{4}t+\ln(t)+\frac{1}{4}  -  \frac {\epsilon}{2}=
$$

\be \label{Re2ndGam}
= \frac{3}{4}\ln\left(  \frac{t}{2 } \right) +\frac{\ln(2\pi)}{2}-\frac{\ln(\pi)}{4}  -\frac{\pi}{4}t+\ln(t)+   \ln \left( \sqrt{ \frac{2 \pi}{t } }  \right) ^\epsilon
\ee

\subsection { Computation of $ \Re_2^-  $}

While  $\Re_2^-$,  the second piece of the Real part using eq. (\ref{sost2}), is:

$$
\Re_2^-  (only \quad till \quad t^{-2} \quad terms)= 
\left( \frac{3-2 \epsilon}{4}  \right)    \left[  \frac{1}{2} \ln \left( 1 +\left( \frac{1+2 \epsilon}{2t} \right)^2 \right)   \right]+
$$

\be\label{SecondoPezzoReGa2}
+\frac{1-2\epsilon}{4} \left[  \quad \frac {2t}{1-2\epsilon}\arctan   \left (         \frac{1-2\epsilon}{2t} \right)    -1\right]
 +\frac{1}{2}\ln\left[    1+ \left(   \frac{1+2 \epsilon}{2t}    \right)^2 \right]=
\ee

$$
\left( \frac{3-2 \epsilon}{4}  \right)    \left[  \frac{1}{2} \left( \left( \frac{1+2 \epsilon}{2t} \right)^2 \right)   \right]+
\frac{1-2\epsilon}{4} \left[  \quad -\frac{1}{3}  \left (         \frac{1-2\epsilon}{2t} \right)^2    \right]
 +\frac{1}{2}\left[     \left(   \frac{1+2 \epsilon}{2t}    \right)^2 \right]+. . . .
$$

Clearly all three terms of $\Re_2^-$ go to zero as $t \rightarrow \infty$.

\subsection { Computation of
$ \Im_1^- $}

For the imaginary part we have:

$$
Im\left( ln\left( e^{-z}z^{   z+\frac{1}{2}   }(2\pi)^{\frac{1}{2}} \pi^{- z}(1-2z) \right )_{z= \frac{1-2 \epsilon}{4} -\frac{i t}{2}} \right)= \Im_1^- + \Im_2^- =
$$

$$
=  -  \left( \frac{t}{2}  \right)    \left[ \ln\left(\frac{t}{2} \right)+ \frac{1}{2} \ln \left( 1 +\left( \frac{1+2 \epsilon}{2t} \right)^2 \right)   \right] +
$$

\be \label{ParteImmag}
+\left( \frac{3-2 \epsilon}{4}  \right) \left[  -\frac{\pi}{2} +\arctan   \left (         \frac{1-2\epsilon}{2t} \right)       \right] +\frac{t}{2}  +\frac{t}{2} \ln(\pi)  +\frac{\pi}{2} -\arctan\left (         \frac{1+2\epsilon}{2t} \right)
\ee

\noindent     Likewise for $t>>\epsilon$ separating real and imaginary parts in a single expression we can write:

$$
\Re_1^- +i \Im_1^- = \frac{3}{4}\ln\left(  \frac{t}{2} \right) + \frac{\ln(2\pi)}{2}-\frac{\ln(\pi)}{4}
$$

\be\label{app8 5}
+ \ln \left( \sqrt{ \frac{2 \pi}{t } }  \right) ^\epsilon -\frac{\pi}{4}t+\ln(t)+ 
 i\left(    - \frac{t}{2} \ln\frac{t}{2\pi}  +\frac{t}{2} +\frac{\pi}{8}+\frac{\pi}{4}\epsilon        \right)
\ee

\subsection { Computation of $\Im_2^-  $}

\be \label{IM2meno}
\Im_2^- (only \quad till \quad t^{-1} \quad  terms) = -      \frac{t}{4} \ln \left( 1 +\left( \frac{1+2 \epsilon}{2t} \right)^2 \right)  +
\left( \frac{3-2 \epsilon}{4}  \right) \left[  \arctan   \left (         \frac{1-2\epsilon}{2t} \right)       \right]  -\arctan\left (         \frac{1+2\epsilon}{2t} \right)=
\ee

$$
= -      \frac{t}{4}  \left( \left( \frac{1+2 \epsilon}{2t} \right)^2 \right)  +
\left( \frac{3-2 \epsilon}{4}  \right) \left[     \left (         \frac{1-2\epsilon}{2t} \right)       \right]  -\left (         \frac{1+2\epsilon}{2t} \right)+. . . .
$$

\subsection { Computation of $\Re_3^- $}

Let us consider the sum in (\ref{Pi&Gamma2}), with substitution  (\ref{sost2}), for K=3.

\noindent    Here too  it is useful to consider the real part:

$$
\Re_3^-  (only \quad till \quad t^{-1} \quad terms)=\Re\left( \sum_{k=1}^{K-1} \frac{B_{2k}   }{  2 k (2k-1) z^{2k-1}} \right)_{z= \frac{1-2 \epsilon}{4} -\frac{i t}{2}  \quad ; \quad K=3 } = 
$$

\be \label{tillB4RE2}
=\frac{1-2 \epsilon}{48\left [        1+ \left( \frac{1-2 \epsilon}{2 t}  \right)^2 \right]    \left( \frac{t}{2} \right) ^2} 
- \frac{  (1-2 \epsilon)^3} {360 t^6  \left(  1+ \left( \frac{1-2 \epsilon}{2 t}  \right)^2   \right)^{3}  }
+\frac{ (1-2 \epsilon)}{30 t^4   \left(  1+ \left( \frac{1-2 \epsilon}{2 t}  \right)^2   \right)^{3}   }=
\ee

$$
=\frac{1-2\epsilon}{12 t^2}-...
$$
\subsection { Computation of $\Im_3^- $}

For the imaginary part:

$$
 Im_3^-  (only \quad till \quad t^{-1} \quad terms) =\Im\left( \sum_{k=1}^{K-1} \frac{B_{2k}   }{  2 k (2k-1) z^{2k-1}} \right)_{z= \frac{1-2 \epsilon}{4} -\frac{i t}{2} \quad ; \quad K=3} = 
$$

\be \label{tillB4IM2}
= + \frac{1}{6 t  \left(  1+\left(   \frac{1-2 \epsilon}{2 t}  \right)^2   \right) }
+ \frac{1}{45 t ^3 \left(  1+\left(   \frac{1-2 \epsilon}{2 t}  \right)^2   \right)^3 }
-\frac{(1-2 \epsilon)^2}{60 t ^5 \left(  1+\left(   \frac{1-2 \epsilon}{2 t}  \right)^2   \right)^3}
\ee

$$
=-\frac{1}{6t}+. . . 
$$

\subsection {Error  merging   } \label  {Error  merging   }

Let us merge the more significant errors (i.e. second and third part ).

   \noindent     For imaginary part we  take only till   $t^{-1}$ term ( for  $\epsilon =0$ we get $\pm \frac{1}{48}$  as in [p. 120 \cite  {Edwards:1974cz}  (1) ]
 
 \be \label {Im2e3Pot1}
 \Im_2^\pm +  \Im_3^\pm  =\frac{1}{t} \left[  \mp \frac{1}{6} \pm \frac{(1+2\epsilon)^2}{16 } \mp \frac{(3\pm2\epsilon)(1 \pm 2\epsilon)}{8} \pm \frac{1 \mp 2\epsilon)}{2}  \right]+. . .  =
 \ee

\noindent      For real part we take only till $t^{-2}$  term:
 
 \be \label {Re2e3Pot2}
 \Re_2^{\pm}+\Re_3^{\pm}  = \frac{1}{t^2}\left[ \frac{3 \pm 2 \epsilon}{4}\frac{1}{2}\left(\frac{1+2\epsilon}{2}\right)^2 -\frac{1\pm 2\epsilon}{12}\left(\frac{1\pm 2\epsilon}{2}\right)^2 +\frac{1}{2}\left(\frac{1\pm 2\epsilon}{2}\right)^2 + \frac{1\pm2\epsilon}{12 }\right] + . . . =
 \ee

 \noindent       The error  (\ref{MaggErrSerAsint1}) contribution, if we limit to   $t^{-1}$ for $\Im$  and $t^{-2}$ for $\Re$, does not appear.

 \noindent       So we have:

 \be \label {ImErr23}
\Im_2^+ +  \Im_3^+=  \frac{ 1 -84 \epsilon +10 \epsilon^2}{ 48 t} + [..] t^{-3} + [..] t^{-5}  + . . .\quad ; \quad
\Im_2^- +  \Im_3^- = - \frac{ 1 +108 \epsilon-12 \epsilon^2}{ 48 t} + [..] t^{-3} + [..] t^{-5}  + . . 
\ee

\noindent      and

 \be \label {ReErr23}
\Re_2^+ +  \Re_3^+=  \frac{ 27 +94 \epsilon +84 \epsilon^2     +8  \epsilon^3  }{ 96 t^2}  + [..] t^{-4} + . . . \quad ; \quad
\Re_2^- +  \Re_3^- =  \frac{ 27 -22 \epsilon+36 \epsilon^2  -8 \epsilon^3    }{ 96 t^2}  + [..] t^{-4}  + . .
\ee

In $\Im_3$ we coud include also the error whose we know  only the bound  ( \ref{RestoInApp} ).

\subsection {Putting pieces together}

\noindent      Let us try to summarize.

\noindent       The  coefficient of first sum of (\ref{Edward(5)pag138})  is:

$$
\left [  \ln( \Gamma(z+1) )= \ln\left( e^{-z}z^{   z+\frac{1}{2}   }(2\pi)^{\frac{1}{2}}  \right )  +\left( \sum_{k=1}^{K-1} \frac{B_{2k}   }{  2 k (2k-1) z^{2k-1}} \right) \right]_{z= \frac{1+2 \epsilon}{4} +\frac{i t}{2}  \quad ; \quad K=3 } =
$$

\be \label{CoeffPrimaSommatoria}
=\Re_1^+ +\Re_2^+ +\Re_3^+ +i( \Im_1^++\Im_2^+ +\Im_3^+ )
\ee

\noindent     while for the second sum and the remainder of  (\ref{Edward(5)pag138}):

$$
 \left [ \ln( \Gamma(z+1) )= \ln\left( e^{-z}z^{   z+\frac{1}{2}   }(2\pi)^{\frac{1}{2}}  \right )  +\left( \sum_{k=1}^{K-1} \frac{B_{2k}   }{  2 k (2k-1) z^{2k-1}} \right) \right]_{z= \frac{1-2 \epsilon}{4} -\frac{i t}{2}  \quad ; \quad K=3 } =
$$

\be \label{CoeffSecSommatoria}
=\Re_1^-+\Re_2^-+\Re_3^- +i( \Im_1^-+\Im_2^-+\Im_3^-)
\ee

\noindent     Let us neglect now $\left( \begin{matrix}  \Re \\ \Im  \end{matrix}\right)_j^+$ and   $\left( \begin{matrix}  \Re \\ \Im  \end{matrix}\right)_j^-$  with  $j>1$  , namely second and third piece.

\noindent    It is useful to isolate real terms that depend only on t.

\noindent  We define:
\be \label {FiDiT}
\ln[F(t)]:=\frac{3}{4}\ln\left(  \frac{t}{2} \right) +\frac{\ln(2\pi)}{2}-\frac{\ln(\pi)}{4}-\frac{\pi}{4}t+\ln(t)
\ee
 $F(t)$ is the scale factor :
 

\be \label {scaleFactorF}
F(t)e^{i\frac{\pi}{4}\epsilon}= e^{\frac{ \ln(2 \pi)}{2}-\frac{ \ln(\pi)+3 \ \ln(2)}{4}} t^{\frac{7}{4}} e^{i \frac{\pi}{4}   (\epsilon+it)} =
 \left(\frac{\pi}{2}\right)^{0.25}t^{\frac{7}{4}} e^{i \frac{\pi}{4}   (\epsilon+it)} 
 \ee


\noindent        In       \cite[p.~176]{Edwards:1974cz}, for  $ \epsilon=0$ and  and  in\cite[p.~5]{AleksandarIvic:2003}

       \be \label {smallScaleFactor}
       f(t)=  \frac{1}{2} \pi^{-1/4}(t^2+ 1/4)|\Gamma(1/4+it/2)| =e^{ Re[ \ln\Gamma(s/2)]}  \pi^{-1/4} \frac{-t^2-1/4}{2}
       \ee 
 
    \noindent   While
    
      $$ F(t)=\Re \left[ (s-1)\Gamma\left( \frac {s}{2}+1\right)  \pi^{-\frac {s}{2}}  \right]_{\epsilon=0}   = f(t)$$
    
    \noindent  So $F(t)$ used here  and    $f(t)$ in   in       \cite[p.~176]{Edwards:1974cz}  and  in\cite[p.~5]{AleksandarIvic:2003} are the same.


\noindent    Note
  that  same  scale factor $F(t)$, applies to (\ref{asintoticaper gamma1re}) (i.e. $\Re_1^+$), as well as to  (\ref{app8 5})(i.e. $\Re_1^-$) .

\noindent      Looking at the imaginary terms and we define :

\noindent For  (\ref{asintoticaper gamma1im})
$$
i \Im_1^+ =
i\left(    \frac{t}{2}\ln\frac{t}{2\pi}  -\frac{t}{2} -\frac{\pi}{8}+\frac{\pi}{4}\epsilon\right)
:=i\left(\theta_1(t)+\frac{\pi}{4}\epsilon  \right)
$$

\noindent    For (\ref{app8 5})
$$
i \Im_1^- =
i\left(   - \frac{t}{2}\ln\frac{t}{2\pi} +\frac{t}{2} +\frac{\pi}{8}+\frac{\pi}{4}\epsilon\right)
:=i\left(-\theta_1(t)+\frac{\pi}{4}\epsilon  \right)
$$

\noindent Where 
\be \label {ThetaDef}
\theta_1(t):=\frac{t}{2}\ln\frac{t}{2\pi}  -\frac{t}{2} -\frac{\pi}{8} = \frac{t}{2} \ln\left(  \frac{t}{2 \pi e}  \right)- \frac{\pi}{8}
\ee

\noindent    Meanwhile referring to    (1)        \cite[p.~120]{Edwards:1974cz}   we have  coherently:
 $$
 \theta(t) =-[  - \theta_1(t) +\Im_2^-(\epsilon = 0) + \Im_3^-(\epsilon = 0) ]=  \theta_1(t) +\Im_2^+(\epsilon = 0) + \Im_3^+(\epsilon = 0)
 $$

\noindent     Looking at the real terms that depend on $\epsilon$ and on t at the same time we have:

\noindent       For  (\ref{asintoticaper gamma1re})
$$
\frac{\epsilon}{2} \left(    \ln  \frac{t}{2}   -\ln(\pi)   \right)=\ln\left[ \left(  \frac{t}{2\pi } \right)   ^{\frac{\epsilon}{2}} \right]
$$

\noindent    For (\ref{app8 5})
$$
\frac{\epsilon}{2} \left(  \ln(\pi )- \ln \left(  \frac{t}{2}   \right)   \right)=
\ln\left[ \left(  \frac{2\pi }{t} \right)   ^{\frac{\epsilon}{2}} \right]
$$

\noindent      Note that also  this apply to \ref{asintoticaper gamma1re} (i.e. $\Re_1^+$), as well as to (\ref{app8 5})  (i.e. $\Re_1^-$). So we could write:

$$
\Re_1^+ = \ln[F(t)] + \ln\left[ \left(  \frac{t}{2 \pi} \right)^{\epsilon / 2}\right]
$$

$$
\Re_1^- = \ln[F(t)] + \ln\left[ \left(  \frac{2 \pi}{t} \right)^{\epsilon / 2}\right]
$$

\noindent       Thus (for $t>>1$) we can use  only the first piece  ($\Re_1$ and $\Im_1$). 

\noindent So we have:

$$
\ln\left[ 
\Gamma \left( \frac{s}{2} +1\right)(1-s)\pi^{-\frac{s}{2}}  \right] 
$$

\be \label{AsintCoeffPrimaSommatoria}
 =\ln\left( e^{-z}z^{   z+\frac{1}{2}   }(2\pi)^{\frac{1}{2}} \pi^{- z}(1-2z) \right )_{z= \frac{1+2 \epsilon}{4} +\frac{i t}{2}} \approx \Re_1^+ +i \Im_1^+
=
 \ln[F(t)] + \ln\left[ \left(  \frac{t}{2 \pi} \right)^{\epsilon / 2}\right] + i\left(\theta(t)+\frac{\pi}{4}\epsilon  \right)
\ee
\noindent      and

$$
\ln\left[  
\Gamma\left( \frac{1-s}{2}+1 \right)(+s)\pi^{-\frac{1-s}{2}}  \right] 
$$
 \be \label {AsintCoeffSecSommatoria}
=\ln\left( e^{-z}z^{   z+\frac{1}{2}   }(2\pi)^{\frac{1}{2}} \pi^{- z}(1-2z) \right )_{z= \frac{1-2 \epsilon}{4} -\frac{i t}{2}} \approx  
\Re_1^-  +i \Im_1^- = 
\ln[F(t)] + \ln\left[ \left(  \frac{2 \pi}{t} \right)^{\epsilon / 2}\right] +i\left(-\theta(t)+\frac{\pi}{4}\epsilon  \right)
\ee

\section { Asymptotic expressions}

Putting together all the terms  in  (\ref{Edward(5)pag138}) and setting $\theta(t) \approx \theta_1(t)$ we have:

$$
-\xi\left( \frac{1}{2}+\epsilon +it \right) \sim  F(t)e^{i\frac{\pi}{4}\epsilon}  
  \sum_{n=1}^{N} \left(   \frac{\sqrt{\frac{t}{2\pi }}}{n} \right)^\epsilon  \frac{e^{i(\theta(t)-t \ \ln(n)}}{\sqrt{n}} 
$$

\be\label {espressione1}
 + F(t)e^{i\frac{\pi}{4}\epsilon}  \left\{  \sum_{n=1}^{N}    \left(\frac {n}{\sqrt{\frac{t}{2\pi }}}  \right)^\epsilon \frac{e^{-i(\theta(t)-t \  \ln(n)}}{\sqrt{n}}  +R(t)  \right\}   \quad ;\quad t>>\epsilon
\ee

\noindent    where $ F(t) =  \left(\pi/2 \right)^{0.25}t^{\frac{7}{4}} e^{- \frac{\pi}{4}  t}  $  , see (\ref  {FiDiT}).

\noindent      Let us look for an exact expression. Using  (\ref{CoeffPrimaSommatoria}) and   (\ref{CoeffSecSommatoria}), and using:

\be \label {Apm}
A^{\pm} = e^{\Re_2^{\pm} + \Re_3^{\pm} + \Re_{error}^{\pm} +i ( \Im_2^{\pm} + \Im_3^{\pm} + \Im_{error}^{\pm}  )  }
\ee

\noindent    we obtain:

$$
-\xi \left ( \frac{ 1}{2}+\epsilon +it \right)  =F(t)e^{i \frac{\pi}{4}  \epsilon}  . .
$$

$$
. . \left\{  \left(   \sqrt{\frac{t}{2\pi }} \right)^\epsilon
  \left\{   
 \sum_{n = 1}^N \left[ \frac{  A^+  n^{-\epsilon} +A^-(2 \pi n/t)^\epsilon     }{\sqrt{n}} \cos(\theta_n) 
 + i  \sin(\theta_n)  \frac{ A^+n^{-\epsilon} -A^-(2 \pi n/t)^\epsilon}{\sqrt{n}}   \right]   \right \} \right\} +
$$

\be \label  {ZSinECosPrecisa} 
. . +F(t)e^{i \frac{\pi}{4}  \epsilon} A^- \left\{ R(t,\epsilon) + \int_{L_0,L_2,L_3} \frac{(-x)^{s}e^{-Nx}dx}{(e^x-1)x}  - \gamma_{L-L_1} \right\}
\ee

\noindent       Where $ \gamma_{L-L_1}$ is evaluated in (\ref  {GammaRef}),
and  $ R(t)$ is
 the   quantity called R in   (6)    \cite[p.~147]{Edwards:1974cz} and
   \cite[p.~154]{Edwards:1974cz} . 
    
    \noindent  Here $ R(t,\epsilon)$ is  the same quantity computed in Appendix \ref{appendix1} also for $\epsilon \ne 0$.

\noindent   The overall bound for  $\int_{L_0,L_2,L_3}$ is shown in Appendix \ref{appendix1}, see (\ref{UpperBounsSum}) and fig. \ref  {MaggiorazioniL023}.

\noindent   Note however that  (\ref  {ZSinECosPrecisa})  cannot be actually computed  because of  the errors in (\ref {Apm})  and in (\ref{UpperBounsSum}), we only have an upper  bound. 

\noindent  Note that  $e^{- i \theta(t)}$, with the same meaning of (1) in   \cite[p.~120]{Edwards:1974cz}      (see also \cite[p.~146]{Edwards:1974cz}  bottom)  is just  included in R; see   in   (6)        \cite[p.~139]{Edwards:1974cz}  .

 \noindent  In       \cite[p.~154]{Edwards:1974cz}],                          the remainder $R(t)$, of Riemann  Siegel Z for  $\epsilon =0$ ,  is computed .

                                                                                                                                                                                                                                                       \noindent  In Appendix \ref{appendix1} the coefficients $C_n(p,\epsilon)$ for $ \epsilon>0$  are generalized.
                                                                                                                                                                                                                                                     
                                                                                                                                                                                                                                                   \be \label {RQuasiTutto}
                                                                                                                                                                                                                                                    R_M(t,\epsilon) \approx (-1)^{N-1} \left( \frac{2 \pi}{t}\right)^{1/4} \sum_{j=0}^M C_j(p,\epsilon)\left( \frac{2 \pi}{t}\right)^{j/2} \quad \quad p= \sqrt{\frac {t}{2 \pi}} - \left \lfloor                \sqrt{\frac {t}{2 \pi } }                  \right \rfloor
                                                                                                                                                                                                                                                    \ee

                                                                                                                                                                                                                                                      \noindent  with:
 \be \label {CZero} 
 C_0(0.5)= 0.382683 .. \le C_0(p)= \frac{(\cos(2 \pi(p^2-p-1/16))}{  \cos(2 \pi p)} \le \cos(\pi/8) \approx 0.923879..=C_0(0)=C_0(1)
 \ee
where $C_0(p,\epsilon) \equiv C_0(p) $ is independent of $\epsilon$; see Appendix \ref{appendix1}                                                                                                                                                                                                                                                                                                                                                                                                                                                                                                                                                                                                                                                                                                                                                                  
                                                                                                                                                                                                                                                                                                                                                                                                                                                                                                                                                                                                                                                                                                                                                                   

\noindent      So we can  put  $ \forall  \epsilon$:
\be \label {R0DiT}
R_0(t)= (-1)^{N-1} \left( \frac{2 \pi}{t}\right)^{1/4}  C_0(p) \quad ; \quad \forall \epsilon
\ee

                                                                                                                                                                                                                                                      \noindent     In the following we use $R_i(t)$ in the meaning of   (\ref  {RQuasiTutto})  with sum in square brackets only  till $C_i$. We use simply $R(t)$  instead of $R_\infty(t)$ .

\noindent   Let us define an useful function  strictly related to $\xi(s)$:

 \be \label {ZRimSiegestesa}
  Z(t,\epsilon):= \sum_{n=1}^{N} \left(   \frac{\sqrt{\frac{t}{2\pi }}}{n} \right)^\epsilon  \frac{e^{i(\theta(t)-t \ln(n)}}{\sqrt{n}} + \sum_{n=1}^{N}    \left(\frac {n}{\sqrt{\frac{t}{2\pi }}}  \right)^\epsilon \frac{e^{-i(\theta(t)-t   \ln(n)}}{\sqrt{n}}  +R(t)
   \sim   \frac{ -  \xi(\frac{1}{2} +\epsilon+it)  } { F(t) e^{i\frac{\pi}{4}\epsilon} } 
\ee

\noindent  By the way: as  $A^{\pm} \rightarrow 1$  ( see (\ref {Apm})),  for  $\epsilon  \rightarrow 0$, then  ( \ref  {ZSinECosPrecisa}) or  (\ref {espressione1} ), lead us to the classical Z of Riemann Siegel   \cite[p.~139]{Edwards:1974cz}.


 
\noindent      Note that $Z(t,\epsilon) e^{i \frac{\pi}{4}  \epsilon}$  has  ``almost ``  the same phase of $ -  \xi(\frac{1}{2} +\epsilon+it) $.

 \noindent     So now it is clear  why we changed the sign in (\ref{Edward(5)pag138})  with respect to (5)  \cite[p.~138]{Edwards:1974cz}.

 
  \noindent  This choice is congruent with $\theta(t)$ in      \cite[p.~119-120]{Edwards:1974cz}                          where  $\theta(t)$ is developped in powers of $1/t$ and $(1/t)^3$. 
  
  \noindent In fact , from ( \ref{ParteImag1} ), we have  
   
   $$
    \theta(t,\epsilon =0)=  \Im\left(   \ln \left[ \Gamma \left(     \frac{s}{2} +1 \right)\pi^{\frac{-s}{2}}(1-s)\right] \right)_{s=1/2+it} =
    \frac{t}{2}\ln\frac{t}{2\pi}  -\frac{t}{2} -\frac{\pi}{8}+ \frac{1}{48t}+\frac{7}{5760t^3} . . .  = \theta(t) 
    $$
 \noindent   like in    (1)  \cite[p.~120]{Edwards:1974cz}]  .

\noindent            An interesting way to express  (\ref   {ZRimSiegestesa}) is using hyperbolic functions:

$$
\frac{-\xi(t,\epsilon)}{F(t)  e^{i\epsilon\frac{\pi}{4}}}  \sim  Z(t,\epsilon)=
2 \sum_{n=1}^N \frac{    \cosh \left[ \epsilon \  \ln\left( \sqrt{\frac {t}{2 \pi n^2}  }\right) \right]   }{\sqrt{n}}
\cos\left(   t \  \ln\left(    \sqrt{\frac {t}{2  e \pi n^2} } \right) - \frac{\pi}{8}   \right) + . .
$$

\be \label {ZSinhECosh}
 . .+ 2 i  \sum_{n=1}^N\frac{  \sinh\left[ \epsilon \  \ln\left( \sqrt{\frac {t}{2 \pi n^2}  }\right) \right]   }{\sqrt{n}}
\sin\left(   t \ \ln\left(    \sqrt{\frac {t}{2  e \pi n^2} } \right) - \frac{\pi}{8}  \right )         +R(t,\epsilon)
\ee

\noindent          Note that, as 
$\frac{\partial}{\partial t} \left[\epsilon \  \ln\left( \sqrt{\frac {t}{2 \pi n^2}  }\right) \right] =\frac{\epsilon }{2 t}   $, each addend of the sum in  (\ref  {ZSinhECosh}) verifies for small  $\frac{\epsilon}{t}$ :
 
 \be \label {RiemanCauchyONE}
 \frac {\partial Im}{ \partial \epsilon }= - \frac {\partial Re}{ \partial t } \quad ;
  \quad
   \frac {\partial Re}{ \partial \epsilon }=  \frac {\partial Im}{ \partial t }
 \ee

 \noindent        For reminder term we have:
 $$
  \frac{\partial R_0(t)}{\partial  t} =\frac{  \partial  \{  (-1)^{N-1}    \left(   \frac{2 \pi}{t}\right)  ^{1/4}   \left[     C_0(p)+. . .\right]  \}}{ \partial t}
  $$

  \noindent       where:
  
  $$
  p=\sqrt{\frac{t}{2 \pi}} -N \quad \rightarrow  \frac{\partial p}{\partial t} =\frac{1}{2\sqrt{2 \pi t} }
  $$
 \noindent       So we have:

   \be \label  {dRsudT_rivisto}
    \frac{\partial R_0(t)}{\partial  t} =
     (-1)^{N-1} 
      \left \{ 
       -  \frac{1}{4t}  \left(   \frac{2 \pi}{t}\right)  ^{1/4} C_0(p)  
         +        \frac{ \partial  C_0(p)}  {\partial p}  \left(   \frac{2 \pi}{t}\right)  ^{1/4}  \frac{1}{2\sqrt{2 \pi t} }  
        \right\}  \approx _{t>>1} . . .
   \ee
$$
 . . .  \approx _{t>>1}   
 (-1)^{N-1} 
      \frac{ \partial  C_0(p)}  {\partial p}  \left(   \frac{2 \pi}{t}\right)  ^{1/4}  \frac{1}{2\sqrt{2 \pi t} }  
$$

\noindent         in  (\ref {Cuno})

\be \label {Cuno2}
 C_1(p,\epsilon)\omega^1 =-\epsilon
  \psi^{(1)}(p) \frac{i \omega}{4 \pi}    
- \psi^{(3)}(p) \frac{\omega}{2^5 \pi^2 3} 
\ee
as $ \omega= \sqrt{\frac {2 \pi}{ t}}$;  see (\ref  {OMEGADef}).

\noindent            As $ \psi^{(1)}(p) =  \frac{ \partial  C_0(p)}  {\partial p}$; see ( \ref  {Psi0}). 
Besides
 $\frac {\partial [\psi^{(3)} \omega ]}{\partial t}
 <<_{t>>1}  \frac{ \partial  C_0(p)}  {\partial p}          \frac{1} {2\sqrt{ 2 \pi t} }  $.

 \noindent          So   taking  only $C_1(p,\epsilon)$ , after $C_0(p)$, in (\ref {RQuasiTutto}):
 
\be \label {R1} 
 R_1(t,\epsilon) = (-1)^{N-1} \left( \frac{2 \pi}{t}\right)^{1/4}  \left[     C_0(p)+C_1(p,\epsilon)\left( \frac{2 \pi}{t}\right)^{1/2}   \right]
 \ee
  
  \noindent      if $t>>1$,       we almost  match,  also for the remainder  $R_1(t,\epsilon)$,   the first  equation of  (\ref {RiemanCauchyONE}) .

\noindent            From  (\ref  {dRsudT_rivisto}) the main part of $ \frac {\partial  Re [R(t,\epsilon)]}{ \partial t} $ is:

$$ \frac {\partial Re [R(t,\epsilon)]}{ \partial t } 
 \approx _{t>>1}
 (-1)^{N-1} 
      \frac{ \partial  C_0(p)}  {\partial p}  \left(   \frac{2 \pi}{t}\right)  ^{1/4}  \frac{1}{2\sqrt{2 \pi t} }$$
While

\be \label {CauchyRiemannForC1}
 \frac {\partial Im [R(t,\epsilon)]}{ \partial \epsilon } 
 \approx _{t>>1}
  -(-1)^{N-1}  \left( \frac{2 \pi}{t}\right)^{1/4}  
  \frac{ \omega}{ 4 \pi}  \frac{ \partial  C_0(p)}  {\partial p} =
  - (-1)^{N-1} 
      \frac{ \partial  C_0(p)}  {\partial p}  \left(   \frac{2 \pi}{t}\right)  ^{1/4}  \frac{1}{2\sqrt{2 \pi t} }
 \ee

\noindent       So   the first of the  (\ref  {RiemanCauchyONE}) is almost verified  also for $R_1(t,\epsilon)$; see (\ref {R1} ).

  \noindent       Note that the second  of ( \ref  {RiemanCauchyONE}) is identically zero for $\epsilon \rightarrow 0$.
 
  \noindent       With reminder $R_1(t,\epsilon)$  which    uses (\ref {RQuasiTutto}) with only $C_0(p) $, and ,$C_1(p,\epsilon)$,
  we can  then  say   that  (\ref  {ZSinhECosh})  , for  $t>>1$ and for 
  small  $\left| \frac{\epsilon}{t} \right|$,  is  almost  holomorphic as the RiemanCauchy (\ref  {RiemanCauchyONE}) holds only  with small discrepancies  that tend to zero as  $t \rightarrow \infty$.

\appendix

\section{ Main integral in \ref{Edward(5)pag138} along path $L_1$ with offset from critical line ($\epsilon > 0$)
 } \label{appendix1}

The aim of this appendix is to  compute the main integral  (\ref  {MainInt}) following the same logical  path  of  \cite[p.~138-155]{Edwards:1974cz}        with the only addition of $\epsilon \ne 0$. 

\noindent     First footnote in \cite[p.~138]{Edwards:1974cz}    suggests exactly this.
A brief summary  shows the path to follow.

\subsection {Summary of logical path }

Let us summarize the logical  path of the following evaluation.

\begin{itemize} 
 \item   To ignore the integral contribution on $L_0,  L_2   ,L_3$;  see (\ref {OnlyL1}).

 \item   To develop $(-x)^{\epsilon-1/2+it}  e^{-Nx}$ as a power series centered at  $a=i \sqrt{2 \pi t}$;  see  (\ref   {XaCentered2}).
 
 \item  To use a ``matching factor''   (\ref{GxMenoASerie}) (rewritten here for easy reading)
 \be \label  {GXmenoA}
g(x-a) = \sum_{n=0}^{n=\infty} b_n (x-a)^n=e^{\left[-\frac{i (x-a)^2}{4 \pi}  -(p+N)(x-a) + (-\frac{1}{2} +\epsilon +it) \ \ln\left( 1 +\frac{x-a}{a}  \right)   \right]  }
\ee

to compute the $b_n(\omega,\epsilon)$ for a perfect matching of all like powers of  $\omega= \sqrt{\frac{2 \pi}{t}}$  obtained by iteration (\ref {IbConEps});   see also (\ref {Sviluppo}).

 \item   To expand the expression obtained by the matching (i.e. $\sum_{n=0}^{n=\infty} b_n (x-a)^n$), on all  $L_1$  line, irrespective of the radius of convergence of the power series centered in $a=i \sqrt{2 \pi t}$  (\ref  {WholeLine}) by a wise contour integration    \cite[p.~147-148]{Edwards:1974cz}   due to Riemann.
 

This allows us  to compute exactly, 
on complex plane $u= x-2 \pi N i$    \cite[p.~147-148]{Edwards:1974cz}, the main  integral  for  $ b_0$ i.e.  $C_0(p)= \psi^{(0)}(p)$  (\ref  {WholeLine}) .
 In order to  compute $C_n(p ,\epsilon)$  by formula  (\ref {ReminderOfXi2}) we need till $b_{3n}$ see  (\ref   {CgrandiCpiccoli}).


\item  through (\ref   {CgrandiCpiccoli}) all pieces are joined to compute $C_i(p,\epsilon)$ of  (\ref {RQuasiTutto}) .

\item  The computation has been actually carried on only for $i=1 ,2$; see (\ref  {Cuno})  and  (\ref   {Cdue}).

  \end{itemize}

\subsection {Main integral over  $L_0+L_1 +L_2+ L_3$}

With (\ref  {ZSinECosPrecisa})  and  (\ref {Apm})  we can write:

$$
 F(t)e^{i \frac{\pi}{4}  \epsilon} A^{-}R(t,\epsilon)=\frac{ (+s)
\Gamma \left(  \frac{1-s}{2} +1\right)\pi^{\frac{-(1-s)}{2}}   }{  (2\pi)^{s-1}  2 \sin(\pi s/2) 2\pi i   }\int_{L_0,L_1,L_2,L_3} \frac{(-1)(-x)^{s-1}e^{-Nx}dx}{e^x-1} 
$$

$$
 =e^{\Re_2^-+\Re_3^- + \Re^{-}_{error} +i( \Im_2^-+\Im_3^-  + \Im^{-}_{error}  ) } . . .
$$
\be \label  {MainInt}
 . .\frac{ F(t)e^{i\frac{\pi}{4}\epsilon} e^{-\pi t/2} e^{-i\theta_1(t)} }{   (  2 \pi  )^{\frac{1}{2} +\epsilon +it} (e^{-i\frac{\pi}{2} \epsilon}-i e^{i\frac{\pi}{2} \epsilon} e^{-t \pi}   ) e^{-i\frac{\pi}{4}}} 
\left ( \sqrt{  \frac{2 \pi }{t}  } \right)^{\epsilon} 
\int_{L_0,L_1,L_2,L_3} \frac{(-x)^{s-1}e^{-Nx}dx}{e^x-1}
\ee

Some note on ( \ref  {MainInt} ).

 \begin{itemize}

\item      Suppose now  (\ref  {ImErr23})  and   (\ref   {ReErr23} ):    $ e^{\Re_2^-+\Re_3^- + \Re^{-}_{error} +i( \Im_2^-+\Im_3^-  + \Im^{-}_{error}  ) }  \approx 1$ .  

\noindent Below in (\ref  {ReminderOfXi2}) for higher order approximation 
$ e^{+i(\Im_2^-+\Im_3^-   ) }  $ is taken in the evaluation of $C_2(p,\epsilon)$;  see (\ref {RQuasiTutto})).

       \item From (\ref {AsintCoeffSecSommatoria})
 $(+s)\Gamma
 \left(  \frac{1-s}{2}+1 \right)\pi^{\frac{-(1-s)}{2}}    \approx  e^{\Re_1^- +i \Im_1^-} =
 F(t)e^{i \frac{\pi}{4}  \epsilon}  e^{-i\theta_1(t)} \left(  \sqrt{  \frac {2 \pi}{t}    } \right)^\epsilon $ .

 \item  The (-1) factor at the integrand numerator stems  for the gathering of $ x$, , under  $(s-1)$ exponent: $ \frac {(-x)^s dx}{(e^x-1)x} =\frac {(-1)(-x)^{s-1} dx}{e^x-1} $.

\item  This $(-1)$ at numerator integrand is, afterward,  joined to $2 i \sin(\pi s/2) $ development  \cite[p.~139]{Edwards:1974cz}         in the denominator  of external coefficient:

  $$
( - 1) 2 i \sin(\pi s/2) = -[e^{i \frac {\pi}{2} \left( \frac{1}{2} + \epsilon +it \right)} -e^{- i \frac {\pi}{2} \left( \frac{1}{2} + \epsilon +it \right)} ]
  =-  e^{- i \frac {\pi}{4} }  [  e^{i \frac {\pi}{2}} e^{  i \frac {\pi  \epsilon }   {2}   } e^{  -\frac {\pi  t }   {2}   }  -    e^{ - i \frac {\pi  \epsilon }   {2}   } e^{  \frac {\pi  t }   {2}   }    ] 
  $$
  
  $$
   = e^{- i \frac {\pi}{4} }  e^{  \frac {\pi  t }   {2}   } [   e^{ - i \frac {\pi  \epsilon }   {2}   } -i e^{  i \frac {\pi  \epsilon }   {2}   } e^{  -\pi  t   }       ]  =
    \frac{ e^{- i \frac {\pi}{4} }  [   e^{ - i \frac {\pi  \epsilon }   {2}   } -i e^{  i \frac {\pi  \epsilon }   {2}   } e^{  -\pi  t   }       ] }{   e^{  \frac {- \pi  t }   {2}   }}
  $$

\item So   in  (\ref  {MainInt}) the  integral is multiplied by: 
$$\frac {1}{(2 \pi)^s} \left( 
F(t)e^{i \frac{\pi}{4}  \epsilon}  e^{-i\theta_1(t)} \left(  \sqrt{  \frac {2 \pi}{t}    } \right)^\epsilon
   \right)
   \left(   
    \frac        {   e^{  \frac {- \pi  t }   {2}   }}  
    { e^{- i \frac {\pi}{4} }  [   e^{ - i \frac {\pi  \epsilon }   {2}   } -i e^{  i \frac {\pi  \epsilon }   {2}   } e^{  -\pi  t   }       ] }  \right)$$

 \end{itemize}

\noindent      If   we put \cite[p.~139]{Edwards:1974cz}:

\be  \label {ExactP0}
P_0 =       \frac{e^{-i \theta_1(t)}  e^{-t \pi/2}     }    {(2 \pi)^{\frac{1}{2}}   (2 \pi)^{ i t }    e^{-i\pi/4}     (1-ie^{-t\pi}  )   }   \quad  for \quad  \epsilon =0 
\ee

\noindent       then ,for $ \forall \epsilon $ we may write: 
\be  \label {ExactPeps}
P_\epsilon =     \left\{   \frac{e^{-i \theta_1(t)}  e^{-t \pi/2}      }    {(2 \pi)^{\frac{1}{2}}  e^{-i\pi/4}    (2 \pi)^{ i t }         (1-ie^{-t\pi}   )  }        \right\}
 \left[    \frac{(1- i e^{-t \pi}  ) }{   (  2 \pi  )^{\epsilon} (e^{-i\frac{\pi}{2} \epsilon}-i e^{i\frac{\pi}{2} \epsilon} e^{-t \pi}   ) } 
   \left ( \sqrt{  \frac{2 \pi }{t}  } \right)^{\epsilon}    
 \right]    \quad  \epsilon  > 0 
\ee

\noindent       As from Appendix \ref{appendix2}  contributions outside $L_1$  are very small  ( see fig. \ref  {MaggiorazioniL023} )  then  (\ref  {MainInt}) can be written: 
\be \label {OnlyL1} 
P_\epsilon \int_{L_0,L_1,L_2,L_3} \frac{(-x)^{s-1}e^{-Nx}dx}{e^x-1} \approx P_\epsilon \int_{L_1} \frac{(-x)^{s-1}e^{-Nx}dx}{e^x-1} 
\ee

\subsection  {Evaluation of  main integral on $L_1$}

We follow the path of          \cite[p.~145]{Edwards:1974cz}   with the only difference of introducing $\forall \epsilon   $.

\noindent         It is useful also           \cite[p.~14-28]{PUGH:1992ar}       .
\noindent    We sum and subtract $a$. 
\be \label  {SumAndSubtr}
-x-a+a = -a \left( 1+ \frac{a-x}{-a}\right) =-a \left( 1+ \frac{x-a}{a}\right)
\ee

\noindent  Where $a= i\sqrt{2 \pi t}$ is the point chosen to expand  (\ref{Numerator}) in powers of $(x-a)$.

\be \label {XaCentered}
(-x)^{-\frac{1}{2} + \epsilon+it} e^{-Nx} =e^{(\epsilon-1/2+it) [ \ \ln(-a)+ \ln\left( 1+\frac{x-a}{a}  \right)} e^{-Na-N(x-a)}
\ee

\noindent       for $\left| \left( \frac{x-a}{a}  \right)\right| <1$

\be \label {LogDevelopment}
 \ln\left( 1+\frac{x-a}{a}  \right)  \approx \frac{x-a}{a} - \frac{(x-a)^2}{2a^2} +\frac{(x-a)^3}{3a^3}- . . . .
 \ee

\noindent                     So:
\be \label {XaCentered2}
(-x)^{-\frac{1}{2} + \epsilon+it} e^{-Nx} =(-a)^{-\frac{1}{2} + \epsilon+it} 
e^{-Na} e^{ \left( \frac{-\frac{1}{2} + \epsilon+it}{a} -N \right)(x-a)  +  (-\frac{1}{2} + \epsilon+it )    \left(   - \frac{(x-a)^2}{2a^2} +\frac{(x-a)^3}{3a^3}- . . . .    \right)    }
\ee

\noindent  We have :
$$
\frac{-\frac{1}{2} + \epsilon+it}{a} = \frac{-\frac{1}{2} + \epsilon+it}{i\sqrt{2 \pi t}} \approx \sqrt{ \frac{t}{2 \pi}   }
$$

\noindent      and:
$$
-\frac{1}{2} \frac{-\frac{1}{2} + \epsilon+it}{a^2} =-\frac{1}{2} \frac{-\frac{1}{2} + \epsilon+it}{( i\sqrt{2 \pi t}    )^2} \approx \frac{i}{4 \pi}
$$

\noindent        The numerator of the integral  can be written  (with  $ a = i \sqrt {2 \pi t}$   \cite[p.~145]{Edwards:1974cz}  bottom):
$$
(-x)^{-\frac{1}{2} + \epsilon+it} e^{-Nx} \approx (-a)^{-\frac{1}{2} + \epsilon+it} e^{-Na} e^{p(x-a)} e^{\frac{i (x-a)^2}{4 \pi}}
$$

\noindent  We look for a factor that transform the above approximate relation in a true equality:

\be \label  {Numerator}
(-x)^{-\frac{1}{2} + \epsilon+it} e^{-Nx} = (-a)^{-\frac{1}{2} + \epsilon+it} e^{-Na} e^{p(x-a)} e^{\frac{i (x-a)^2}{4 \pi}}g(x-a) =
\ee

$$
 =(-a)^{-\frac{1}{2} + \epsilon+it} e^{-Na} e^{p(x-a)} e^{\frac{i (x-a)^2}{4 \pi}}\sum_{n=0}^\infty b_n(x-a)^n
$$

\noindent      Where $a= i\sqrt{2 \pi t}$ is the point chosen to expand ( \ref{Numerator} ) in powers of $(x-a)$.

$g(x-a)$ defined in         \cite[p.~145]{Edwards:1974cz}     bottom and in  \cite[p.~24]{PUGH:1992ar}   is:

\be \label {GxMenoASerie}
g(x-a) = \sum_{n=0}^{n=\infty} b_n (x-a)^n=e^{\left[-\frac{i (x-a)^2}{4 \pi}  -(p+N)(x-a) + (-\frac{1}{2} +\epsilon +it)  \ \ln\left( 1 +\frac{x-a}{a}  \right)   \right]  }
\ee

\noindent      It is apparent that  $b_0 =1 $ even when $\epsilon \ne 0$.

\noindent Because (\ref {GxMenoASerie}) for $x = a$ is bound to be $\epsilon$ independent  as $\epsilon$ effects are multiplied by $ \ln(1)$.

\noindent      For a given $\epsilon$, with  the power series   (\ref{GXmenoA})``$a$'' centered the reminder $R(t,\epsilon)$  is:

$$
R(t,\epsilon)=       P_\epsilon
\int_{L1 } \frac{(-a)^{-\frac{1}{2} + \epsilon+it} e^{-Na} e^{p(x-a)} e^{\frac{i (x-a)^2}{4 \pi}}g(x-a) dx }{e^x-1}  =
$$

\be \label {RtEps}
=   P_\epsilon
(-a)^\epsilon
\int_{L1 } \frac{(-a)^{-\frac{1}{2} +it} e^{-Na} e^{p(x-a)} e^{\frac{i (x-a)^2}{4 \pi}}g(x-a) dx }{e^x-1}  = . .
\ee

\noindent       We have eliminated $F(t) e^{i \frac{\pi}{4}  \epsilon} $  as we are looking for  (\ref  {ZRimSiegestesa})  i.e.$
  Z(t,\epsilon) \sim   \frac{ -  \xi(\frac{1}{2} +\epsilon+it)  } { F(t) e^{i\frac{\pi}{4}\epsilon} } $.

\noindent      If we ignore  $e^{-t \pi}$   with respect to 1  ( note that for  $t>14$  : $e^{-t \pi} < e^{-14 \times \pi} \approx 10^{-19}$) we get:

\be  \label {ApproxP0}
P_0 \approx        \frac{e^{-i \theta_1(t)}  e^{-t \pi/2}     }    {(2 \pi)^{\frac{1}{2}}   (2 \pi)^{ i t }    e^{-i\pi/4}     (1-i\times 0  )   } \quad  for \quad  \epsilon =0 
\ee

\noindent      and for same reason  we have:
\be  \label {ApproxPeps}
P_\epsilon \approx     \left\{   \frac{e^{-i \theta_1(t)}  e^{-t \pi/2}      }    {(2 \pi)^{\frac{1}{2}}  e^{-i \pi/4}    (2 \pi)^{ i t }         (1-i \times 0   )  }        \right\}
 \left[    \frac{(1- i  \times 0   ) }{   (  2 \pi  )^{\epsilon} (e^{-i\frac{\pi}{2} \epsilon}-i  \times 0    ) } 
   \left ( \sqrt{  \frac{2 \pi }{t}  } \right)^{\epsilon}    
 \right]    \quad  for \quad  t >>  \epsilon  > 0 
\ee

\noindent        So $R(t,\epsilon)$ in ( \ref  {RtEps}) becomes:

\be \label {RconEps}
 = P_0  \quad
 \int_{L1 } \frac{(-a)^{-\frac{1}{2} +it} e^{-Na} e^{p(x-a)} e^{\frac{i (x-a)^2}{4 \pi}}g(x-a) dx }{e^x-1}  \quad t>>\epsilon  >0
\ee

\noindent      because:

$$
 P_\epsilon
(-a)^\epsilon =
 P_0 \frac{ 1  }{   (  2 \pi  )^{\epsilon} e^{-i\frac{\pi}{2} \epsilon}} \left ( \sqrt{  \frac{2 \pi }{t}  } \right)^{\epsilon} (-a)^\epsilon=
\left ( \sqrt{  \frac{ 2 \pi t}{2 \pi t}  } \right)^{\epsilon}  e^{i\frac{\pi}{2} \epsilon} (-i)^\epsilon
 P_0  \quad
= P_0  \quad t >>\epsilon >0
$$

 \noindent      What happens is that main integral in  (\ref {MainInt} ) gives rise to a factor  $(-a)^\epsilon= (-i \sqrt{t 2 \pi})^\epsilon$ [  (2) p. 146  \cite{Edwards:1974cz}  ] who keeps invariant the coefficient  $P_0$ before integral   $\forall \epsilon$.

\noindent       We can fetch out from integral in ( \ref  {RconEps} ) the expression :

$$
[-a]^{-\frac{1}{2} + i t} e^{Na} =[ -i \sqrt{2 \pi t}  ] ^{-\frac{1}{2} + i t} e^{-N i \sqrt{2 \pi t}}=
[ -i \sqrt{2 \pi t}  ] ^{-\frac{1}{2} + i t}  e^{-N i 2 \pi \sqrt{\frac{ t}{2 \pi}}}
$$

$$
 =[ -i \sqrt{2 \pi t}  ] ^{-\frac{1}{2} + i t}  e^{-N i 2 \pi \sqrt{\frac{ t}{2 \pi}}} =
  [ -i \sqrt{2 \pi t}  ] ^{-\frac{1}{2} + i t}   e^{-N i 2 \pi   (N+p)}
$$

\noindent      So ( \ref  {RconEps} ) becomes:

\be \label {RconEps2}
 = P_0  \quad (-a)^{-\frac{1}{2} +it} e^{-Na}  
 \int_{L1 } \frac{   e^{p(x-a)} e^{\frac{i (x-a)^2}{4 \pi}}g(x-a) dx }{e^x-1}  \quad; \quad  t>>\epsilon
\ee

\noindent      Developing the coefficient 
taking into account (\ref  {ApproxPeps}) gives the same  result  of \cite[p.~146]{Edwards:1974cz}: 

$$
 P_0  \quad (-a)^{-\frac{1}{2} +it} e^{-Na}= 
 \frac{e^{-i \theta_1(t)}  e^{-t \pi/2}    [ -i \sqrt{2 \pi t}  ] ^{-\frac{1}{2} + i t}   e^{-N i 2 \pi   (N+p)}  }    {(2 \pi)^{\frac{1}{2}}   (2 \pi)^{ i t }    e^{-i\pi/4}      }  
 $$
 
 $$
  = \frac{e^{-i \theta_1(t)} ( e^{i \pi/2} )^{it}   \left [  \sqrt{\frac{t} {2 \pi }} \right ] ^{-\frac{1}{2} + i t} (-i)^{-\frac{1}{2}+it}  e^{-N i 2 \pi   (N+p)}  }    {  (2 \pi)    ( e^{i\pi/2}  ) ^{- \frac{1}{2}}   } =  
  \frac{i  e^{-i \theta_1(t)}   \left (  \sqrt{\frac{t} {2 \pi }} \right ) ^{-\frac{1}{2} + i t}    e^{-N i 2 \pi   (N+p)} }{2 \pi}
 $$
  $$
  = \frac{-  e^{-i \theta_1(t)}   \left (  \sqrt{\frac{t} {2 \pi }} \right ) ^{-\frac{1}{2} + i t}    e^{-N i 2 \pi   (N+p)} }{2 \pi i}
 $$

  \noindent      Which is the  coefficient of  $  \int_{L1 } \frac{ e^{p(x-a)} e^{\frac{i (x-a)^2}{4 \pi}}g(x-a) dx }{e^x-1}   $ in     \cite[p.~146]{Edwards:1974cz}     (apart $1-ie^{-\pi t}$ at denominator  just reduced to 1).
  
 \noindent       We  can change variable  accordingly to    \cite[p.~146]{Edwards:1974cz}:
 
 \be \label {XChangeU}
 x= u+2 \pi i N, \quad x-a = u+2 \pi i N - 2 \pi i \sqrt{\frac{t}{2 \pi}} =u- 2 \pi i p
 \ee
  
 \noindent      So numerator integrand change in  : 
 \be \label {NumChange}
 e^{p(x-a)} e^{\frac{i (x-a)^2}{4 \pi}} \rightarrow e^{p u - 2 \pi i p^2} e^{i \frac{u^2-4 \pi^2 p^2-4 \pi iu p}{4 \pi}} =
 e^{i \frac {u^2}{4 \pi} +2pu -2 \pi i p^2-\pi i p^2 } 
 \ee
  
 \noindent      We can now take away from integral the objects:
  $-2 \pi i p^2-\pi i p^2$
 
 \noindent      Elaborating we get: 
    
  $$
-N i 2 \pi   (N+p)   = - 2 \pi i N^2 - 2 \pi i N p \rightarrow - 2 \pi i N^2 - 2 \pi i N p -2 \pi i p^2-\pi i p^2
= 
  $$
  
\noindent       which is the same as :
 
  $$
 - i \pi (N+p)^2 -i \pi N^2  - 2 \pi i p^2= -i \pi N^2  -i \pi p^2 -i \pi 2 Np -i \pi N^2  - 2 \pi i p^2
  $$

   \noindent      Besides  :$- i \pi (N+p)^2 -i \pi N^2  - 2 \pi i p^2 =  -i \pi N^2-i \pi \frac {t}{2 \pi}     -2 \pi i p^2 $
 
  \noindent      but:
   
  $$
 - \theta(t) =-\left[   \frac {t}{2} ln\left(    \frac{t}{2 \pi e} \right) - \frac {\pi}{8} \right]
  $$
so we have:  
 
 $$
 e^ { -i \left[  
  \frac {t}{2}  \ ln \left(    \frac{t} {2 \pi e} 
  \right)  \right] }
   \left (     \sqrt{   \frac{t} {2 \pi } } \right ) ^ { + i t} e^{ -i \pi \frac {t}{2 \pi} } = e^{i\frac{\pi}{8}}
   $$
  
  $$
  -e^{-i \pi N^2}=(-1)^{N-1}
  $$
 
 \noindent      So like in  (2)     \cite[p.~147]{Edwards:1974cz}  we  can write   $R(t,\epsilon)$ 
      in (  \ref   {RconEps} )  as:
 
 \be \label {LikeEdw4Pg147}
 R(t,\epsilon)= 
  (-1)^{N-1} \left (     \sqrt{   \frac{t} {2 \pi } } \right ) ^ { -\frac{1}{2}}
 \frac{e^{i\frac{\pi}{8}} e^{-2 \pi i p^2}  }{2 \pi i} \int_{L_1} \frac{e^{i \frac {u^2}{4 \pi} +2pu  } g(u-2 \pi i p ) }{e^u-1}du
 \ee
 
  \noindent    The usual notation is:
 
 \be \label {WholeLine}
C_0(p)= \psi(p) =  \frac{e^{i\frac{\pi}{8}} e^{-2 \pi i p^2}  }{2 \pi i} \int_{whole \ \ L_1  \ \ line} \frac{e^{i \frac {u^2}{4 \pi} +2pu  } }{e^u-1}du
 \ee
 
 \noindent     Details of  contour integration used to get  (\ref  {WholeLine}) are in    \cite[p.~147-148]{Edwards:1974cz}.

\subsection {Dependence from $\epsilon$}

So from (5)  \cite[p.~145]{Edwards:1974cz}   and  (3)  \cite[p.~150]{Edwards:1974cz} it  is obvious that  $C_0(p)$,  defined        in      \cite[p.~154]{Edwards:1974cz}is $\epsilon$ independent; see  (\ref{RQuasiTutto}) .

 \noindent   For   $C_i(p,\epsilon) $  with $i>0$ the usefulness of following manipulation is apparent when logarithmic derivative of $g(x-a) $ with exponential  and with the sum: $\sum_{n=0}^{n=\infty} b_n (x-a)^n$ are equated; see ( \ref {IbConEps})


 \noindent      Putting

\be \label {OMEGADef}
\omega = \sqrt{\frac{2 \pi}{t}}
\ee
  \noindent     for \cite[p.~148]{Edwards:1974cz}  (bottom page) we have:

$$
\frac{\partial }{\partial x} \left[ -\frac{i(x-a)^2}{4 \pi } - (p+N)(x-a) + \left(   - \frac{1}{2} +\epsilon + i t \right) \ln\left( 1+\frac{x-a}{a} \right) \right] 
$$

\be \label{Sviluppo}
= \frac{x-a}{2 \pi i} -\omega^{-1} + \frac{\epsilon - \frac{1}{2} +2 \pi i \omega^{-2}}{2 \pi i \omega^{-1} + (x-a)} =
 \frac{(x-a)^2 +2 \pi i ( \epsilon - \frac{1}{2}  )}{ 2 \pi i [ 2 \pi i \omega^{-1} + (x-a)  ]   } = \frac{\sum_{n=0}^{n=\infty} b_n n  (x-a)^{n-1}}{ \sum_{n=0}^{n=\infty} b_n (x-a)^n}
\ee

 \noindent      where, equating like powers of $(x-a)$, it is easy to get  $ b_k$  recursively.

 \noindent     So we get the generalization with $\epsilon >0$ of [(5)    \cite[p.~152]{Edwards:1974cz}:   

\be \label {IbConEps}
b_{n+1} = \omega \frac{2 \pi i (n+ 0.5-\epsilon)b_n -b_{n-2}}{4 \pi^2 (n+1)} 
\quad \rightarrow
 b_n = \left( \frac{i}{2 \pi} \right)\frac{\omega}{n} P_\epsilon(n,n) b_{n-1 }+\left( \frac{i}{2 \pi} \right)^2\frac{\omega}{n} b_{n-3 }
\ee

 \noindent   where $P_\epsilon(n1,n2)$ is defined as:
 \be \label {PEps}
 P_\epsilon(n_1,n_2)= \prod_{n=n1}^{n=n2}\left(   n- \frac{1}{2}  -\epsilon \right)
 \ee

   \noindent    with  

 \noindent   $b_0 =1 \quad ; b_{-1} =  b_{-2} =0$


¨
 \noindent    For higher order approximation  the expression  (3) \cite[p.~150]{Edwards:1974cz}    is exactly:
$$
 [ F(t)e^{i\frac{ \pi}{4}\epsilon    }  ]^{-1}  \quad \frac{ (+s)\Gamma\left(  \frac{1-s}{2}+1 \right)\pi^{ - \frac{1-s}{2}}   }{  (2\pi)^{s-1}  2 \sin(\pi s/2) 2\pi i     }\int_{C_{\delta,N}} \frac{(-x)^{s-1}e^{-Nx}dx}{e^x-1}
$$

\be \label  {ReminderOfXi2}
  =
 R(t,\epsilon)\approx (-1)^{N-1} \left(    \sqrt{\frac{2 \pi}{t}}  \right)^{0.25} e^{i( \Im_2^- + \Im_3^-)}\left[\sum_{k=0}^{k=K} b_k(\omega) c_k \right ]
\ee

 \noindent      where $b_k(\omega,\epsilon) $ is a polynomial of max degree n in $\omega$ ,while $c_k $ (`$\epsilon$ independent ``) is a linear combination of even, if $n$  is even, or odd, if $n$  is odd, derivatives with respect to p of 
$C_0(p)= \psi(p )$; see (\ref  {WholeLine}).

 \noindent   And where the recovered  $\Im_2^- + \Im_3^-$  terms in $t^{-1}$ 
is given in (\ref{Im2e3Pot1}). 
While $c_k$ is from (\ref{Icn}) and $ F(t)$ is defined in (\ref{FiDiT} ).
 
  \noindent With variable change  $x-a = u-2 \pi ip $; see( \ref  {XChangeU} ) and  (\ref  {Numerator}) .

  \noindent  The $c_n$ are  in  (3)             \cite[p.~150]{Edwards:1974cz}   which is the same of   (\ref  {ReminderOfXi2}) here.
 
 \noindent    They are computed in (\ref {Icn})    which derives  from     (4) \cite[p.~150]{Edwards:1974cz} given below :
 
 $$e^{2 \pi i y^2} \sum_{m=0}^\infty\frac{\psi^m(p)}{m!}y^m = \sum_{n=0} ^\infty \frac{(2 y)^n}{n!}  c_n
  \quad (4) \quad \quad \cite[p.~150]{Edwards:1974cz} $$ 
 
  \noindent       equating  homogeneous coefficients in $y$ of left and right series.

 \noindent     Result is the following:

\be \label {Icn}
c_n =  \frac {e^{i \pi/8}  e^{-2\pi p^2} } {2 \pi i} \int_{\Gamma} \frac{e^{i u^2/4\pi}   e^{2pu}}{e^u-1} (u-2 \pi i p)^n du= \frac{n!}{2^n} \sum_{j=0}^{   \left  \lfloor   n / 2  \right \rfloor  }  \frac{   (2 \pi i)^j}       { j!   }    \frac{ \psi^{  (n-2j)   } (p)  }{(n-2j)!   } 
\ee


 \be \label {Icn2}
 = \sum_{j=0}^{   \left  \lfloor  n / 2  \right \rfloor   } c_{n,n-2j} \psi^{(n-2j)}=c_n
\ee

 \noindent    So $c_{k,j}$  is the coefficient of  contribution at $c_k$ of derivative $ \psi^{(j)}$.

  \noindent   Remind that 
\be \label {Psi0}
\psi^{ ( 0)   } (p) = \frac{cos[ 2 \pi(p^2-p-1/16)]}{cos(2 \pi p)}= C_0(p)
\ee

  \noindent    From iteration  (\ref{IbConEps}) we have:

$$ b_{-2}(\omega)= 0 \quad, \quad b_{-1}(\omega)= 0 \quad, \quad b_0(\omega) =1 \quad, \quad b_1(\omega) = B_{1,1 } \omega^1 $$ , 

$$ b_2(\omega) = B_{2,2 } \omega^2 \quad, \quad  b_3 = B_{3,3} \omega^3 +B_{3,1} \omega^1 \quad,
$$
$$ \quad  b_4(\omega) = B_{4,4} \omega^4 +B_{4,2} \omega^2   \quad,
$$
$$
 \quad b_5(\omega) = B_{5,5} \omega^5 +B_{5,3} \omega^3$$

$$ b_6(\omega) = B_{6,6 } \omega^6 + B_{6,4 } \omega^4 + B_{6,2 } \omega^2  \quad, \quad b_7 (\omega)= B_{7,7} \omega^7 +B_{7,5} \omega^5   + B_{7,3} \omega^3   \quad, 
$$
$$
\quad b_8(\omega) = B_{8,8} \omega^8 +B_{8,6} \omega^6   + B_{8,4} \omega^4 
$$ 

$$b_9(\omega) = B_{9,9} \omega^9 +B_{9,7} \omega^7   + B_{9,5} \omega^5  + B_{9,3} \omega^3$$

  \noindent   Where $B_{m,k}(\epsilon)$ is the coefficient of $\omega^k$ in the polynomial $b_m(\omega,\epsilon)$.  They contain the dependence from $\epsilon$  throughout ( \ref  {PEps}) factors.

  \noindent   The  polynomial $b_m(\omega,\epsilon)$ max degree is $n$.

  \noindent   While 
 the minimum  degree is: $$ \left  \lfloor   \frac {n}{3 }  \right \rfloor     +n- 3 \times \left  \lfloor   \frac {n}{3 }   \right \rfloor$$ .


  \noindent   The formula connecting  big single indexed $C_i(p,\epsilon)$ in  (\ref {RQuasiTutto}), small single indexed $c_k$, from ( \ref {Icn}), 
small doubly indexed c : $c_{k,j}$, from  (\ref {Icn2}) and 
doubly indexed  big B above  i.e.   $B_{m,k}(\epsilon)$ is :


\be \label   {CgrandiCpiccoli}
C_k \omega^k =\sum_{k \le i =3k-2j \le 3k \quad j= 0,1,2...} c_i B_{i,k} \omega^k=\sum_{k \le i =3k-2h \le 3k \quad h= 0,1,2...}
 \left[  \sum_{j=0}^{  \left  \lfloor   n / 2 \right \rfloor    } c_{i,i-2j} \psi^{(i-2j)} 
 \right ] B_{i,k} \omega^k
\ee


   \noindent  Once fixed $k$ of $C_k \omega^k $ in  (\ref {RQuasiTutto}  ) , the  index $i$  can  varies between $3k$ and  $3k-2h\ge k$  , with $h= 0,1,2...$.

\subsection {Evaluation of $C_1(p,\epsilon)$ }

From (\ref{PEps} ) we have:

$$
B_{1,1} = \frac{i \omega}{2 \pi}  \frac{P_\epsilon(1,1)}{ 1!}
$$

$$
B_{2,2} =\left( \frac{i \omega}{2 \pi}\right)^2  \frac{P_\epsilon(1,2)}{ 2!}
$$

$$
B_{3,3} =\left( \frac{i \omega}{2 \pi}\right)^3  \frac{P_\epsilon(1,3)}{ 3!}  \quad, \quad B_{3,1} = \left(\frac{ \omega}{2 \pi}  \right)    \frac{i^2}{2 \pi 3}
$$

$$
C_1 \omega^1 = c_1 B_{1,1} \omega^1 + c_3B_{3,1} \omega^1
$$
   \noindent   where
 $$
  c_1 =  \frac{1!}{2^1} \frac{(2 \pi i)^0}{0! } \frac{\psi^{1-0} } {0!}
 $$
and
$$
c_3 = \frac{3 !}{2^3}  \left[   \frac{(2 \pi i)^0}{0! } \frac{\psi^{3-0} } {3!}   +  \frac{(2 \pi i)^1}{1! } \frac{\psi^{3-2} } {1!}   \right] =c_{3,3}\psi^{(3)}+c_{3,1}\psi^{(1)}
$$
  \noindent    Terms like  $\omega^1$ have ``b'' produced by  $b_1$and by  $b_3$, so:

$$
C_1(p,\epsilon) \omega^1=
\frac{i \omega}{2 \pi} P_\epsilon(1,1) \frac{1!}{2^1} \frac{(2 \pi i)^0}{0! } \frac{\psi^{1-0} } {0!}
+
\frac{3 !}{2^3}  \left[   \frac{(2 \pi i)^0}{0! } \frac{\psi^{3-0} } {3!}   +  \frac{(2 \pi i)^1}{1! } \frac{\psi^{3-2} } {1!}   \right]
 \left[  - \frac{\omega}{2 \pi} \quad  \frac{1}{2 \pi 3}  \right] 
$$
\be \label {Cuno}
 = \psi^{(1)}(p) \frac{i \omega}{4 \pi}\left(  \left(\frac{1}{2}-\epsilon\right) - \frac {3!}{2^2}\quad  \frac{2 \pi}{2 \pi 3} \right) - \psi^{(3)}(p) \frac{\omega}{2^5 \pi^2 3} = -\epsilon
  \psi^{(1)}(p) \frac{i \omega}{4 \pi}    
- \psi^{(3)}(p) \frac{\omega}{2^5 \pi^2 3} 
\ee

  \noindent Putting $\epsilon =0$, then
the imaginary part is zero and the real part happens to be like in       \cite[p.~153-154]{Edwards:1974cz}   .

\subsection {Evaluation of $C_2(p,\epsilon)$ }

Consider now:

$$
C_2 \omega^2 = c_2B_{2,2} \omega^2 +c_4 B_{4,2}\omega^2 +c_6B_{6,2} \omega^2 
$$

  \noindent    Terms like $\omega^2$  stem from $b_{3\times 2}$ till  $b_{ 2} $: $b_2(\omega)$, $b_4(\omega)$ e  $b_6(\omega)$,  so may  be involved only the even  derivatives:

$$
\psi^{(m)}(p) \quad; \quad  m= 0,2,4,6
$$

  \noindent     $ B_{2,2}$ is computed above.

$$
B_{4,2} = \left(\frac{ \omega}{2 \pi}  \right)^2    \frac{(i)^{4-1}}{ \pi 4!}[ P_\epsilon(4,4) +  \frac{3	\times 2}{2} P_\epsilon(1,1)]
$$

  \noindent    and

 $$
 B_{6,2} = \left( \frac{ \omega}{2 \pi}\right)^2       \frac{ (i)^{6-2}}   {\pi^2  6!} \frac{5\times 4}{2}
 $$

  \noindent     Besides 

$$
c_2 =\frac{2!}{2^2}\left[ \frac{\psi^{(2)}}{2!}   +\frac{  \psi^{(0)}   }{0!} \frac{(2 \pi i)^1}{1!}\right]=c_{2,2}\psi^{(2)}+c_{2,0}\psi^{(0)}
$$

  \noindent    with:
$$
c_4 =\frac{4!}{2^4}\left[ \frac{\psi^{(4)}}{4!}  +\frac{  \psi^{(2)}   }{2!} \frac{(2 \pi i)^1}{1!} +\frac{  \psi^{(0)}   }{0!} \frac{(2 \pi i)^2}{2!}\right]=c_{4,4}\psi^{(4)}+c_{4,2}\psi^{(2)}+c_{4,0}\psi^{(0)}
$$
  \noindent      and:
$$
c_6 =\frac{6!}{2^6}\left[ \frac{\psi^{(6)}}{6!}  
+\frac{  \psi^{(4)}   }{4!} \frac{(2 \pi i)^1}{1!} 
+\frac{  \psi^{(2)}   }{2!} \frac{(2 \pi i)^2}{2!} +\frac{  \psi^{(0)}   }{0!} \frac{(2 \pi i)^3}{3!}\right]=
c_{6,6}\psi^{(6)}+c_{6,4}\psi^{(4)}+c_{6,2}\psi^{(2)}+c_{6,0}\psi^{(0)}
$$

  \noindent    Let us compute separately terms like  $\omega^2$ produced by $\psi^{(0)}(p),\quad \psi^{(2)}(p),\quad \psi^{(4)}(p),\quad \psi^{(6)}(p)$:


\subsubsection { Factors of $\psi^{(0)}(p)$}

$$
\psi^{(0)}(p)\left\{  \frac{2!}{2^2}   \frac{(2 \pi i)^1}{1!}  \left( \frac{i \omega}{2 \pi}\right)^2  \frac{P_\epsilon(1,2)}{ 2!}    \right\}+
$$

\be \label{PsiZero}
+ \psi^{(0)}(p)\left\{ 
\frac{4!}{2^4}\frac{(2 \pi i)^2}{2!}  \left(\frac{ \omega}{2 \pi}  \right)^2    \frac{(i)^{4-1}}{ \pi 4!} \left[ P_\epsilon(4,4) +  \frac{3	\times 2}{2} P_\epsilon(1,1)\right]   +
\frac{6!}{2^6}\frac{(2 \pi i)^3}{3!}   \left( \frac{ \omega}{2 \pi}\right)^2       \frac{ (i)^{6-2}}   {\pi^2  6!} \frac{5\times 4}{2}           \right\}=
\ee


$$
=\frac{ \psi^{(0)}(p)  \pi \omega^2 i^3}{(2 \pi)^2}
\left\{ 
\frac{P_\epsilon(1,2)}{2} -\frac{1}{2^3}[P_\epsilon(4,4) + 3P_\epsilon(1,1) ] +\frac{5 \times  2^5}{2^8 \times3}
\right\}=
\frac{ \psi^{(0)}(p)   \omega^2 i}{2^5 \pi 3} 
\left\{ 
1+12 \epsilon+. . . 
   \right\}
$$

  \noindent dropping   powers of $\epsilon$ greater that 1

  \noindent  Consider  now the phase shift due the multiplication by $A^-$  of the reminder  integral in \ref  {ZSinECosPrecisa}  .

  \noindent  Let us limit to the  $t^{-1}$ term, then  from ( \ref{Im2e3Pot1} ) we have:

$$
\theta^+(t,\epsilon) = \theta_1(t)- \frac{ -1 +84 \epsilon+10 \epsilon^2}{ 48 t} \quad ; \quad
\theta^-(t,\epsilon) = -\theta_1(t)- \frac{ 1 +108 \epsilon-12 \epsilon^2}{ 48 t} 
$$

  \noindent  from (4)       in     \cite[p.~147]{Edwards:1974cz}, we can argue that $C_n(p,\epsilon)$ terms are   \noindent  multiplied by:

\be \label{CorrNegFase}
e^{ i( \Im_2^-+\Im_3^-)} =
e^{ i\left(    - \frac{ 1 +108 \epsilon-12 \epsilon^2}{ 48 t }   + . . .   \right)} \approx  1- \frac{ i \omega^2(1 +108 \epsilon-12 \epsilon^2    )}{ 2^5 \pi 3}
\ee

   \noindent   If we multiply main integral in  (\ref{MainInt})  by (\ref{CorrNegFase}), then the term  in $\omega^0$ : $\psi^{(0)}(p)$, will be translated to  $\omega^2$ summing up  with homogeneous term  stemming by other source.

  \noindent      In this case( sum of (\ref{PsiZero}) with (\ref {CorrNegFase})  at $\epsilon=0$) we have a cancellation ,so for  $\epsilon=0$ ,$\psi^{(0)}(p)$ does not contribute to  $C_2(p,\epsilon=0)$.


\subsubsection { Factors of $\psi^{(2)}(p)$}

$$
\psi^{(2)}(p)\left\{   \frac{2!}{2^2 2!}  \left( \frac{i \omega}{2 \pi}\right)^2  \frac{P_\epsilon(1,2)}{ 2!}       +  
  \frac{4!}{2^4} \frac {(2 \pi i)^2}{2! 2!} \left(\frac{ \omega}{2 \pi}  \right)^2    \frac{(i)^{4-1}}{ \pi 4!} \left [ P_\epsilon(4,4) +  \frac{3	\times 2}{2} P_\epsilon(1,1)\right] \right\}+
$$
\be \label{Psi2}
+\psi^{(2)}(p)\left\{  
 \frac{6!}{2^6 2! } \frac{(2 \pi i)^2}{2!}\left( \frac{ \omega}{2 \pi}\right)^2       \frac{ (i)^{6-2}}   {\pi^2  6!} \frac{5\times 4}{2}
   \right\}
\ee

$$
 =\frac{\psi^{(2)}(p)   5}{2^6 \pi^2}\left(  1-\frac{ 8 \epsilon}{10} \right) -\frac{  \psi^{(2)}(p)  3   }{2^7 \pi^2} \left[      \left(   1-\frac{2 \epsilon}{3} \right)     \left(  1-2 \epsilon \right)    \right] 
 -\frac{\psi^{(2)}(p)   5 }{ 2^7 \pi^2} 
$$

  \noindent      which is $ \frac{\psi^{(2)}(p)   }{2^6 \pi^2}   $ for $\epsilon=0$.


\subsubsection { Factors of $\psi^{(4)}(p)$}

$$
\psi^{(4)}(p)\left\{   
 \frac{4!}{2^4 4!}   \left(   \frac{ \omega}{2 \pi}  \right)^2    \frac{(i)^{4-1}}{ \pi 4!} \left [ P_\epsilon(4,4) + 
  \frac{3	\times 2}{2} P_\epsilon(1,1)\right]   + 
  \frac{6!}{2^6 4!}\frac{(2 \pi i)^1}{1!}  \left( \frac{ \omega}{2 \pi}\right)^2       \frac{ (i)^{6-2}}   {\pi^2  6!} \frac{5\times 4}{2}
  \right\} 
$$

$$
 =\frac{ -i \psi^{(4)} (p) (5-8 \epsilon)  }{2^9 \pi^3 3} + \frac{ i  \psi^{(4)} (p) 5 }{2^9 \pi^3 3} 
$$
  \noindent    which is 0 for $\epsilon=0$. So $\psi^{(4)}(p)$ gives no contribution to $C_2(p,\epsilon=0)$.


\subsubsection { Factors of $\psi^{(6)}(p)$}

$$
\frac{ 6!}{2^6} \left[  .... \frac{(2 \pi i)^0}{0!}   \frac{\psi^{(6)}(p)}{ 6 !}     ..... \right]  \left( \frac{
 \omega}{2 \pi}\right)^2       \frac{2 \quad 5 (i)^{6-2}}   {\pi^2  6!} = \frac{  \psi^{(6)}(p) 5 \omega^2   }{2^7 \pi^4 6!} =\frac{  \psi^{(6)}(p) 5    }{2^7 \pi^4 6!}
$$

  \noindent     So it does not depend on $\epsilon $ accordingly with third row in \cite[p.~153]{Edwards:1974cz} bottom .

  \noindent      In summary for  $\epsilon =0$ the contribution to  $\omega^2$ term simplify drastically  
  \cite[p.~154]{Edwards:1974cz} bottom:

\be \label {Cdue}
C_2(p,\epsilon =0) =   \frac{\psi^{(2)}(p)   }{2^6 \pi^2}     +\frac{  \psi^{(6)}(p) 5    }{2^7 \pi^4 6!}
\ee

  \noindent     Similarly we can built  all the $C_i(p,\epsilon)$ of    (\ref {RQuasiTutto}) :

\be \label{Cdue}
C_2(p,\epsilon)=\frac{  i  \psi^{(0)}(p)     }{ 2^5 \pi 3} \left[  15  \left(  1-\frac{ 2 \epsilon}{3}  \right)(1-2 \epsilon) - 5 -9\left( 1-\frac{8 \epsilon}{10}  \right)      \right] - \frac{ i  \psi^{(0)}(p)   (1 +108 \epsilon-12 \epsilon^2    )}{ 2^5 \pi 3}  + 
\ee

$$
 \frac{\psi^{(2)}(p)   5}{2^6 \pi^2}\left(  1-\frac{ 8 \epsilon}{10} \right) -\frac{  \psi^{(2)}(p)  3   }{2^7 \pi^2} \left[      \left(   1-\frac{2 \epsilon}{3} \right)     \left(  1-2 \epsilon \right)    \right] 
 -\frac{\psi^{(2)}(p)   5 }{ 2^7 \pi^2} + \frac{ -i \psi^{(4)} (p) (5-8 \epsilon)  }{2^9 \pi^3 3} + \frac{ i  \psi^{(4)} (p) 5 }{2^9 \pi^3 3}  +\frac{  \psi^{(6)}(p) 5    }{2^7 \pi^4 6!}
$$

\clearpage

\section{Bounds for Main Integral  on paths:  $L_0,L_1,L_2, L_3      $. With offset from critical line ($\epsilon > 0$)} \label{appendix2}

\subsection {Upperbound of Main integral on $L_0$ path}

We follow  \cite[p.~141-142]{Edwards:1974cz}  `` Estimation of the integral away from saddle point'', first case  j=0, with the only difference of introducing $\epsilon > 0$.

  \noindent    We move  in the complex plane x  on line  : $x = a+ k e^{i\frac{\pi}{4}}$
with real $ k \ge1$, and , we choose as boundary between  $L_0$ and  $L_1$ point:
\be \label {BoundaryL0L1}
Boundary(L_0,L_1) = a+ \frac{|a|}{r} e^{i\frac{\pi}{4}} \quad ; \quad   1<r \le 2
\ee

  \noindent     Where  $a$ is  defined in ( \ref {a}).

  \noindent     In \cite[p.~141-142]{Edwards:1974cz}  $ r = 2$, but,  if $r>1$ 
it is enough for the logarithmic series   (\ref {LogDevelopment})    to converge.

  \noindent  If we write the numerator module  of the integrand in  (\ref{MainInt} )  as  $  | (-x)^{s-1}e^{-Nx}  | =e^{\phi(k)}$ where
 $$
 \phi(k) =Re[ (\epsilon-1/2+it) \  ln(a+ k e^{i\frac{\pi}{4}})-N(a+ k e^{i\frac{\pi}{4}})]
$$

  \noindent    then
$$
\frac{d \phi(k)  }{d k} =Re\left[   N  e^{i\frac{\pi}{4}}\left( \frac{-\frac{1}{2} +\epsilon +it}{N( a+ k e^{i\frac{\pi}{4}}  )  }    -1\right) \right]=
Re\left[   N  e^{i\frac{\pi}{4}}\left( \left[ \frac{-\frac{1}{2} +\epsilon }{N}+\frac{it}{N} \right]\frac{1}{ i\sqrt{2\pi t}+ \frac{k}{\sqrt{2}} +i \frac{k}{\sqrt{2}}  }  -1  \right)   \right]=
$$

$$
Re\left[   N  e^{i\frac{\pi}{4}}\left( \left[ \frac{-\frac{1}{2} +\epsilon }{N}+\frac{it}{N} \right]
\frac{- i\sqrt{2\pi t}+ \frac{k}{\sqrt{2}} -i \frac{k}{\sqrt{2}} }{ (\sqrt{2\pi t}+ \frac{k}{\sqrt{2}} )^2+\frac{k^2}{2}  }    -1   \right)    \right]=
$$

$$
=\frac{N}{\sqrt{2}} \left[  Re\left\{  \left[ \frac{-\frac{1}{2} +\epsilon -t }{N}+i\frac{t -\frac{1}{2} +\epsilon}{N} \right]\frac{ \frac{k}{\sqrt{2}}- i\sqrt{2\pi t} -i \frac{k}{\sqrt{2}} }{ (\sqrt{2\pi t}+ \frac{k}{\sqrt{2}} )^2+\frac{k^2}{2}  } \right\}   -1 \right]
$$

$$
=\frac{N}{\sqrt{2}} \left[   \frac{ \left ( -\frac{1}{2} +\epsilon -t \right)  \frac{k}{\sqrt{2}}     }
{
N\left( (\sqrt{2\pi t}+ \frac{k}{\sqrt{2}} )^2+\frac{k^2}{2}   \right)}   + \frac{   \left(   t -\frac{1}{2} +\epsilon  \right)    \left(  \sqrt{2\pi t} + \frac{k}{\sqrt{2}}   \right)  }
{  N\left( (\sqrt{2\pi t}+ \frac{k}{\sqrt{2}} )^2+\frac{k^2}{2}   \right)    }  -1  \right]  = . .
$$

\be \label {PendNeg}
=\frac{N}{\sqrt{2}} \left[   \frac{ \left ( -\frac{1}{2} +\epsilon  \right)  \frac{ 2 k}{\sqrt{2}}   + \left(   t -\frac{1}{2} +\epsilon  \right)    \left(  \sqrt{2\pi t}   \right) }
{
N\left( (\sqrt{2\pi t}+ \frac{k}{\sqrt{2}} )^2+\frac{k^2}{2}   \right)  }    -1  \right] < 0
\ee
Where   to be $ . . <0$ is justified by:

   (1) $N= \sqrt{\frac{t}{2 \pi}}$ ,

    (2)   Both  numerator  and denominator  have : $  \sqrt{\frac{t}{2 \pi}}( 2\pi t)$.
  
    \noindent  At numerator, for $|\epsilon| <1/2$ and  $ k \ge1$ we have to add negative quantities, while at denominator we add all positive quantities .

   \noindent  In \cite[p.~141-142]{Edwards:1974cz} it  is chosen $r=2$, but , in  order  (\ref  {LogDevelopment} ) to converge    $ r >1  $  is enough.


  \noindent    So  the numerator module  of the integrand in (\ref{MainInt})  on $L_0$, is  max in $ x= Boundary(L_0,L_1) = a+ \frac{|a|}{r} e^{i\frac{\pi}{4}} $; see  [\ref   {BoundaryL0L1} .

  \noindent    As $ \ln(w) = \ln|w| + i\angle [w]$ and  $|e^w|= e^{Re[w]}$ ,with complex $w$,
So

$$
 e^{  \left(  \epsilon-\frac{1}{2}   \right)  \ln\left|   -a-\frac{|a|}{r} e^{i\frac{\pi}{4}} \right |} =  
 \left| \sqrt{2 \pi t}  \right|   ^{      \epsilon-\frac{1}{2}     }
  \left|   \left( -i \left(    1 + \frac{1}{r \sqrt{2}}  \right) - \frac{1}{r \sqrt{2}} \right)    \right|  
$$

  \noindent    then  the numerator module  of the integrand in  \ref{MainInt} is at most:

\be \label {NumModuleIntegrand}
e^{    \Re\left[     \left(  \epsilon-\frac{1}{2}  +it \right) \
\ln \left(  -a-\frac{|a|}{r} e^{i\frac{\pi}{4}} \right) 
  -N \left( a+\frac{|a|}{r} e^{i\frac{\pi}{4}} \right)     \right]   } =
  e^{  \left(  \epsilon-\frac{1}{2}   \right)  \  \ln\left|   -a-\frac{|a|}{r} e^{i\frac{\pi}{4}} \right| -t \Im \left  [ \ ln\left(  -a-\frac{|a|}{r} e^{i\frac{\pi}{4}} \right)  \right]-N \frac{|a|}{r\sqrt{2}}        } = . . .
\ee
  
    \noindent      remembering $a= i\sqrt{2 \pi t}$  and  $ N+p =  \sqrt{ \frac{t}{2 \pi}}  $, we have: 
  
  $$-N \frac{|a|}{r\sqrt{2}}  =   - \left(\sqrt{ \frac{t}{2 \pi}} -p \right) \frac{\sqrt{2\pi t}}{r\sqrt{2}}     = -\frac {t}{r\sqrt{2}}  + p \frac{\sqrt{\pi t}}{r} \quad \quad \quad  0\le p <1 $$ 
 
    \noindent     Or in other words:
    we can elaborate $-N\frac{\sqrt{\pi t}}{r}$ to get :

$$
-N\frac{\sqrt{\pi t}}{r} < -(N-1)\frac{\sqrt{\pi t}}{r} \le  -\left(\sqrt{\frac{t}{2 \pi}}-1\right)\frac{\sqrt{\pi t}}{r} =  - \frac{t}{r \sqrt{2}} + \frac{\sqrt{\pi t}}{r}
$$ 
 
  \noindent    besides 

$$
 \Im \left[  \   \ln\left(  -a-\frac{|a|}{r} e^{i\frac{\pi}{4}} \right) \right]  = \Im \left[ \  \ln(a)   \right]   +
  \Im\left[  -1+ \frac{i}{r} e^{i\frac{\pi}{4}}  \right]  
  =
   -\frac{\pi}{2}  +
 \Im\left[  -1+ \frac{i}{r} \left(   \frac {1}{\sqrt{2}} +  \frac {i}{\sqrt{2}}  \right)  \right]  
$$

  \noindent       so we can write:

$$
e^{   -t \Im \left[ \  \ln\left(  -a-\frac{|a|}{r} e^{i\frac{\pi}{4}} \right) \right] -N \frac{|a|}{r\sqrt{2}}                            }
= 
e^{-t\left(   -\frac{\pi}{2} -\arctan\left(   \frac{   \frac{1}{r \sqrt{2}} }{\frac{1}{r \sqrt{2}}  +1}     \right)\right) -N\frac{\sqrt{\pi t}}{r} } 
$$
  
    \noindent      defining further :
 \be \label {K0}
 K_0(r) :=\arctan\left( \frac{1}{r\sqrt{2}+1} \right) -\frac{1}{r\sqrt{2}}   \quad ; 1< r \le 2
\ee
  
  \noindent       we can write ( \ref {NumModuleIntegrand}) as:
  
$$
 =\left| \sqrt{2 \pi t}  \right|   ^{ \epsilon-\frac{1}{2} } \left|   -i\left(  \left(    1 + \frac{1}{r \sqrt{2}}  \right) -i \frac{1}{r \sqrt{2}} \right)    \right|  
e^{-t\left(   -\frac{\pi}{2} -\arctan\left(   \frac{   \frac{1}{r \sqrt{2}} }{\frac{1}{r \sqrt{2}}  +1}     \right)\right) -N\frac{\sqrt{\pi t}}{r} }  
$$
  \noindent       taking $ p= 1$, we can bound  the numerator module  of the integrand in  (\ref{MainInt}  )i.e.(\ref {NumModuleIntegrand}  ) with:
\be \label {Ineq1}
 <   \left| \sqrt{2 \pi t}  \right|   ^{ \epsilon-\frac{1}{2} } 2e^{t\frac{\pi}{2}} e^{t \left( \arctan\left( \frac{1}{r\sqrt{2}+1} \right) -\frac{1}{r\sqrt{2}}     \right) } e^{\frac{\sqrt{\pi t}}{r}} 
 =2
 \left| \sqrt{2 \pi t}  \right|   ^{ \epsilon-\frac{1}{2} } 
  e^{t\frac{\pi}{2}}
 e^{t K_0(r)}     e^{\frac{\sqrt{\pi t}}{r}} \quad ; 1< r \le 2
\ee

  \noindent      because 

$$e^{-N\frac{\sqrt{\pi t}}{r}} <e^{ - \frac{t}{r \sqrt{2}} + \frac{\sqrt{\pi t}}{r}}$$

  \noindent     but:

$$ 1<\left|     \left(    1 + \frac{1}{r \sqrt{2}}  \right) -i \frac{1}{r \sqrt{2}}    \right| < 2  \quad \quad \forall r : 1<r \le 2 $$

  \noindent        We can bound, by ( \ref  {Ineq1} ),   the numerator  value on $L_0$ with the value  that $|(-x)^{s-1}e^{-Nx}|$ assumes  on $ Boundary( L_0,L_1)$, so the integral on $L_0$  is bounded  by:

   \noindent       $u= \Re(x)$ from  $u= \frac{ |a|} { r \sqrt{2}} =\frac{\sqrt{\pi t}}  {r}  $   to $\infty$ 
  \be \label {DaCmplxAReal}
  \int_{L0}^\infty \frac{|dx|}{|e^x-1 |} \le \int_{\frac{\sqrt{t \pi}}{r}}^\infty \frac {e^{-u }  \sqrt{2 } du  }   {1-e^{-u}} \quad; \quad  u = \Re(x)\quad; \quad  x= i \sqrt{2 t \pi} +\frac{ \sqrt{2 t \pi} } {r}e^{i \pi/4} k\quad; \quad |dx|= \sqrt{2}du
  \ee
 because :
  $ |1-e^{-x}| =  \sqrt{ (1-e^{-\Re(x)} cos (-Im(x)) ^2 +(e^{-\Re(x)} sin (-Im(x) )^2} \ge | 1-e^{-\Re(x)} cos (-Im(x)|   \ge | 1-e^{-\Re(x)} | = 1-e^{-u}$  . And we e have  $\Re(x) >0$

   \noindent      In other words we  change an integral on ``x''  complex variable  with an integral on real variable $u$  which gives us a bound of module.
 
   \noindent    And we can write
 
 \be \label {BoundOnInt}
  \int_{\frac{\sqrt{\pi t}}{r}}^\infty  \frac {|dx|}{|e^x-1|} \le   \int_{\frac{\sqrt{\pi t}}{r}}^\infty  \frac {e^{-u} \sqrt{2}du}{1-e^{-u}}
   <
    \frac{           \sqrt{2}   [-e^{-u}]_{    \frac{   \sqrt{\pi t}   }  {r}  }    ^{\infty}            }
 {    1-e^{      -\frac{\sqrt{\pi t}}  {r}        }     } =
   \frac{          \sqrt{2}    e^{- \frac{   \sqrt{\pi t}   }  {r} }         }
 {    1-e^{      -\frac{\sqrt{\pi t}}  {r}        }     } 
  \ee

    \noindent    Afterword $ \int_{L0}^\infty \frac{|dx|}{|e^x-1 |}$ is   multiplied by   \ref  {Ineq1}.

   \noindent      Then the  $ A^-R(t,\epsilon) e^{i\frac{\pi}{4}\epsilon} $ in   \ref{MainInt}  module  (   divided by $F(t)$ for \ref  {ZRimSiegestesa}   ) is bounded above  ( because of    ( \ref  {PendNeg})    ) by:
$$
\left|e^{\Re_2^-+\Re_3^- +i( \Im_2^-+\Im_3^-)}\frac{
e^{i\frac{\pi}{4}\epsilon} 
e^{-\pi t/2} e^{-i\theta_1(t)} }{   (  2 \pi  )^{\frac{1}{2} +\epsilon +it} |e^{-i\frac{\pi}{2} \epsilon}-i e^{i\frac{\pi}{2} \epsilon} e^{-t \pi}   |  } 
\left ( \sqrt{  \frac{2 \pi }{t}  } \right)^{\epsilon} \right|\left| \sqrt{2 \pi t}  \right|   ^{ \epsilon-\frac{1}{2} } 
2
  e^{t\frac{\pi}{2}}   
    e^{K_0(r)t} 
      e^{\frac{\sqrt{\pi t}}{r}}        \int_{L_0} \frac {|dx|}{|e^x-1|} 
$$

$$
 =
 \frac{e^{(Re_2^-+Re_3^- )} }{   (  2 \pi  )^{\frac{1}{2} +\epsilon } |e^{-i\frac{\pi}{2} \epsilon}-i e^{i\frac{\pi}{2} \epsilon} e^{-t \pi}   | } 
\left ( \sqrt{  \frac{2 \pi }{t}  } \right)^{\epsilon}    \left| \sqrt{2 \pi t}  \right|   ^{ \epsilon-\frac{1}{2} } 
 2
 e^{K_0(r) t}  e^{\frac{\sqrt{\pi t}}{r}} 
   \int_{L_0} \frac {|dx|}{|e^x-1|} 
$$

$=\left[ \sqrt{2 \pi } \right]^\epsilon$  simplifies with  $\frac {1}{ (  2 \pi  )^{\frac{1}{2} +\epsilon } } $ resulting in :
$ \frac {1}{ (  2 \pi  )^{\frac{1+\epsilon}{2}  } } $ 

$$
=
 \frac{e^{(Re_2^-+Re_3^- )} }{   (  2 \pi  )^{\frac{1+\epsilon}{2}  } |e^{-i\frac{\pi}{2} \epsilon}-i e^{i\frac{\pi}{2} \epsilon} e^{-t \pi}   |      } 
  \left| \sqrt{2 \pi t}  \right|   ^{ -\frac{1}{2} } 
2
 e^{K_0(r)t}  e^{\frac{\sqrt{\pi t}}{r}} 
   \int_{L_0} \frac {|dx|}{|e^x-1|}\approx   . .
    $$
    
   $$
    \quad t>10 \quad  , \quad \forall \epsilon  \quad \quad  |e^{-i\frac{\pi}{2} \epsilon}-i e^{i\frac{\pi}{2} \epsilon} e^{-t \pi}   |   \rightarrow 1
$$

$$
 . .\approx   \frac{e^{(Re_2^-+Re_3^- )} }{   (  2 \pi  )^{\frac{1+\epsilon}{2}  }     } 
  \left| \sqrt{2 \pi t}  \right|   ^{ -\frac{1}{2} } 
  2
 e^{K_0(r)t}  e^{\frac{\sqrt{\pi t}}{r}} 
   \int_{\frac{\sqrt{\pi t}}{r}}^\infty  \frac {|dx|}{|e^x-1|}  < . . .
$$

  \noindent     See (\ref  {ReErr23}  )
  and (\ref  {ImErr23}) ) 
  $\Re_2^-+\Re_3^- = \left(  \frac{27}{96} +f(\epsilon) \right) \frac{1}{t^2} + \left(  ...  \right)\frac{1}{t^4} ...$ so they can be neglected , for $ t >10$  and applying \ref  {BoundOnInt} we have 

$$
 . .<   \frac{    \left| \sqrt{2 \pi t}  \right|   ^{ -\frac{1}{2} }   }{   (  2 \pi  )^{\frac{1+\epsilon}{2}  }     } 
  2
 e^{K_0(r)t}  e^{\frac{\sqrt{\pi t}}{r}} 
   \frac{          \sqrt{2}    e^{- \frac{   \sqrt{\pi t}   }  {r} }         }
 {    1-e^{      -\frac{\sqrt{\pi t}}  {r}        }     } 
$$

  \noindent      So we reach an  stricter  upper bound  with respect to \cite[p.~143]{Edwards:1974cz}; see  (\ref {MainInt} ). So we have:

$$
 \left|   \frac{ (+s)\Gamma\left(  \frac{1-s}{2} +1\right)\pi^{\frac{-(1-s)}{2}}   }{  (2\pi)^{s-1}  2sin(\pi s/2) 2\pi i   }\int_{L_0} \frac{(-x)^{s-1}e^{-Nx}dx}{e^x-1} \right|  < . . .
$$

\be \label   {MaggL0}
  . . .  <
    \frac{    \left| \sqrt{2 \pi t}  \right|   ^{ -\frac{1}{2} }   }{   (  2 \pi  )^{\frac{1+\epsilon}{2}  }     } 
  2
 e^{K_0(r)t}  e^{\frac{\sqrt{\pi t}}{r}} 
   \frac{          \sqrt{2}    e^{- \frac{   \sqrt{\pi t}   }  {r} }         }
 {    1-e^{      -\frac{\sqrt{\pi t}}  {r}        }     } <
 \frac { 2 \sqrt{2}  e^{K_0(r) t}}{    \left[\sqrt{2 \pi } \right]^{ (1+\epsilon )/2 }     \sqrt{ \sqrt { 2 \pi t   } } }
 = UpperBound_{L_0}(t)
\ee

  \noindent   Note that  upper bound on $L_0$  (\ref {MaggL0} ) is  decreasing  while  $\epsilon$ increasing .

\subsection { $L_1$ spurious contribution} \label {Spurious}

To compute ( \ref {RconEps2})  in Appendix \ref{appendix1}  we used, like   \cite[p.~147]{Edwards:1974cz}, a contour integral  which gives an exact result (see (\ref  {CZero}) and (\ref   {LikeEdw4Pg147})  ) on entire line $ L \supset L_1$ .

  \noindent     We are interested in the evaluation between : (\ref {BoundaryL0L1}) $
Boundary(L_0,L_1) = a+ \frac{|a|}{r} e^{i\frac{\pi}{4}} \quad ; \quad   1<r \le 2$  and ( \ref {BoundaryL1L2} ): $Boundary(L_1,L_2) = a- \frac{|a|}{r} e^{i\frac{\pi}{4}} \quad ; \quad   1<r \le 2$, where $a$ is the saddle point given in (  \ref {a} ). 

  \noindent     Here we evaluate the spurious  contribution inserted with this contour integral in $L-L_1$,  that in  ( \ref {ZSinECosPrecisa} ), is referred as :
\be \label {GammaRef}
 \gamma_{L-L_1} 
\ee

   \noindent     For standard Gaussian formulas see \cite[p.~183]{MurrayRSpiegel:2004cz}    .

$$ 
G(x)=\left( \frac{c}{\pi} \right)^{1/2}e^{-cx^2}    = \frac {1}{\sqrt{2 \pi} \sigma} e^{-\frac{x^2}{2 \sigma^2}}  \quad  ;\quad c=\frac {1}{2 \sigma^2}
 $$
  \noindent     An approximate expression for large values of argument ( $x> 2 \sigma$ )   can be obtained from : 
 \be \label {ErfStuff}
  \int_{-x}^{+x} G(x) dx =Erf \left(  y=  \frac{x}{\sigma \sqrt{2}} \right) \approx 1-\frac {e^{-y^2}}{\sqrt{\pi} y}
\left( 
1-\frac{1}{2y^2}+ \frac {1\times 3}{(2y^2)^2}- \frac {1 \times 3 \times 5}{ (2y^2)^3 }+ . . .
  \right) \quad for \quad  y >\sqrt{2}
 \ee

   \noindent     Now let us focus  (\ref {RconEps2}) on : 
\be \label {RiemanContourInt}
  \int_{L1  \rightarrow L} \frac{ e^{p(x-a)} e^{\frac{i (x-a)^2}{4 \pi}}g(x-a) dx }{e^x-1}   
  \ee
 \noindent     See also      \cite[p.~146]{Edwards:1974cz}      .

 \noindent    putting:

\be \label {YDef}
x-a =\sqrt {i} y = y\left(       \frac{1}{ \sqrt{2}} + \frac{i}{ \sqrt{2}}          \right)
\ee

 \noindent  namely:

\be \label {XFunY
}x= a +\sqrt {i} y = i \sqrt{2\pi t} +\sqrt {i} y =  2\left( 
 \frac{y}{2 \sqrt{2}} + i \left( \sqrt{\frac{\pi t} {2}} + \frac{y}{2 \sqrt{2}}   \right) \right)
\ee

 \noindent     As  $(\sqrt{i})^2=i$ we have:

$$
e^{\frac{i (x-a)^2}{4 \pi}} =e^{\frac{i (\sqrt{i} y)^2}{4 \pi}}=e^{-\frac{ ( y)^2}{4 \pi}}
$$

 \noindent    and: 
$$
e^{p(x-a)}=e^{p \left (   \frac{y}{ \sqrt{2}}  + i \frac{y}{ \sqrt{2}}       \right)}
$$

 \noindent      Let us pose $ g(x-a) =1$ and compute:

$$
 \left|   \frac{ e^{p(x-a)} e^{\frac{i (x-a)^2}{4 \pi}}}{e^x-1} \right|
 $$

 \noindent      for  $ \frac{1}{e^x-1} $ on the line  $L_1$ we have

$$
 \frac{ 1  }{   e^{x}-1     }  = 
  \frac{     e^{ \left( -\frac{1}{2}\right) x }} { e^{x/2}-e^{-x/    2}  } 
  =
2 \frac
{
e^{ 
\left(-1\right)
 \left(  \frac{y}{2 \sqrt{2}  } + i \left( \sqrt{\frac{\pi t} {2}} + \frac{y}{2 \sqrt{2}}   \right)   \right) }        }
 {  \sinh \left( \frac{y}{2 \sqrt{2}} + i \left( \sqrt{\frac{\pi t} {2}} + \frac{y}{2 \sqrt{2}}   \right) \right)}  
$$

 \noindent     While:
$$
 \sinh \left(\frac{y}{2 \sqrt{2}} + i \left( \sqrt{\frac{\pi t} {2}} + \frac{y}{2 \sqrt{2}}   \right) \right)
 $$
 $$
 = \sinh \left( \frac{y}{2 \sqrt{2}}\right) \cos \left(   \sqrt{\frac{\pi t} {2}} + \frac{y}{2 \sqrt{2}}    \right) +
 i \cosh \left( \frac{y}{2 \sqrt{2}}\right) \sin \left(   \sqrt{\frac{\pi t} {2}} + \frac{y}{2 \sqrt{2}}    \right) 
$$

 \noindent     and taking the module:
$$
\left|2 \frac
{
 e^{ 
 \left( - \frac{y}{2 \sqrt{2}  } - i \left( \sqrt{\frac{\pi t} {2}} + \frac{y}{2 \sqrt{2}}   \right)   \right) }        }
 {  \sinh \left( \frac{y}{2 \sqrt{2}} + i \left( \sqrt{\frac{\pi t} {2}} + \frac{y}{2 \sqrt{2}}   \right) \right)} \right|  \le
  2 \frac {
  e^{ 
 \left( - \frac{y}{2 \sqrt{2}  }    \right)
  }
   } 
     {  \sqrt{
      \sinh^2 \left( \frac{y}{2 \sqrt{2}}\right) \cos^2 \left(   \sqrt{\frac{\pi t} {2}} + \frac{y}{2 \sqrt{2}}    \right) +
  \cosh ^2\left( \frac{y}{2 \sqrt{2}}\right) \sin^2 \left(   \sqrt{\frac{\pi t} {2}} + \frac{y}{2 \sqrt{2}}    \right)        }
      }
 $$

 \noindent     If we avoid $y=0$ and the points $\in L_1$, as :

  $\sinh[\alpha] < \cosh[\alpha] \quad \forall  \alpha \in \Re$,

\be \label {ExpressioneFratta}
\left|2 \frac
{
 e^{ 
 \left( - \frac{y}{2 \sqrt{2}  } - i \left( \sqrt{\frac{\pi t} {2}} + \frac{y}{2 \sqrt{2}}   \right)   \right) }        }
 {  \sinh \left( \frac{y}{2 \sqrt{2}} + i \left( \sqrt{\frac{\pi t} {2}} + \frac{y}{2 \sqrt{2}}   \right) \right)} \right |_{ y \not\in L_1} 
 <\left [
  2 \ \frac{   e^{ 
 \left( - \frac{y}{2 \sqrt{2}  }    \right)
 } }
 { \left |\sinh \left( \frac{y}{2 \sqrt{2}}\right) \right |} 
 \right]_{ y \not\in L_1}   \left\{
  \begin{matrix}  
      < 2 \quad  &   for  \quad  y < 0 \quad \\
 \quad \quad  \quad  \approx  2  e^{ 
 \left(     - \frac  {y}   { \sqrt{2}  }    \right)
  }  \quad &for \quad   y \quad > 0  \\
 \end{matrix}  
 \right\}
 \ee

 \noindent     For others factors of integrand,by (\ref {YDef} ):

\be \label {ModGau}
 \left|    e^{p(x-a)} e^{\frac{i (x-a)^2}{4 \pi}} \right| = \left|    e^{py\left(       \frac{1}{ \sqrt{2}} + \frac{i}{ \sqrt{2}}\right) } e^{\frac{-y^2}{4 \pi}} \right|  = 
 e^{
 \left(           \frac {py} { \sqrt{2}} -\frac{ y^2}{4 \pi}     \right)
  } 
 \ee

 \noindent  The module ( \ref {ModGau}) has a Gaussian shape, with  $\sigma = \sqrt{2 \pi}\approx  2.5066 . .   $, and, the max at 
 \noindent   $ y_{max} = \frac {2 p \pi}{\sqrt{2}} \approx p \times4.4428 . . = p  \sigma \times 1.772 . . \quad ; \quad 0 \le p \le 1$

 \noindent   While ( $r \approx 1$) : $-|a| <y < |a|= \sqrt{2\pi t} \approx \sqrt{t} \times 2.5066 .            . = \sqrt{t} \times \sigma$

 \noindent   Or, if we  follow    \cite[p.~141]{Edwards:1974cz}      i.e. $r=2$: $-\frac{|a|}{2} <y < \frac{|a|}{2}$ .

 \noindent     If  $ t =20 $,  $p = \sqrt{t/(2 \pi)}-  
  \left \lfloor               \sqrt{t/(2 \pi) }               \right \rfloor=0.7841 . .$

 \noindent     So we have the maximum of Gaussian shape  at:
  $$p  \sigma \times 1.772 . .=  \sigma \times (0.7841 . \times  1.772 .)= 1.3894 \times  \sigma   $$

 \noindent     See (\ref {a}):  $|a| = \sqrt{2\pi t} \approx \sqrt{t} \times  \sigma =  4.472 \sigma$.

 \noindent     So for $r\ge 1$ we have a distance of Gaussian maximum to  the boundary  given by:
$$
(4.472-1.38) \sigma= 3.19 \sigma
$$
¨
 \noindent      So let us call
 $$Contribution_{y \not \in L_1} = \int_{y \not \in L_1} e^{
 \left(           \frac {py} { \sqrt{2}} -\frac{ y^2}{4 \pi}     \right)
  } dy $$
  \noindent       the  additional spurious contribution  that enters in contour integration of (\ref  {RiemanContourInt}), both for positive or negative $y$, taking into account (\ref   {ExpressioneFratta}).

  \noindent         And lets call   
  $$Contribution_{y  \in L_1} = \int_{y  \in L_1} e^{
 \left(           \frac {py} { \sqrt{2}} -\frac{ y^2}{4 \pi}     \right)
  } dy $$ the contribution of the integral in $-|a| <y < |a|$ for $ t =20$. 

  \noindent     Then we have :  
 \be \label {ContributionInOut}
 |Contribution_{y  \in L_1} | \times 0.003 >| Contribution_{y \not \in L_1} |  
 \ee
 
 \noindent  because
   (\ref {RiemanContourInt}), Gaussian beyond $3 \sigma$  . 
If we take $r=2 $ like in        \cite[p.~141]{Edwards:1974cz}       the boundary is at $\frac{|a|}{2} =2.237 \sigma$ and we ought to choose a $t> 20$  in order to reach the same $3\sigma $ distance of border from Gaussian maximum .

  \noindent   Anyhow   with $t>100$ (\ref  {ContributionInOut}) holds  with huge margin .

  \noindent  This analysis together with  fig. \ref {MaggiorazioniL023}
  and (\ref    {UpperBounsSum}), explains why  the case study of  $\xi(t,\epsilon=0) $, in           \cite[p.~155]{Edwards:1974cz}        presents a very good match with Haselgrove tabel       \cite[p.~122]{Edwards:1974cz} at least   from  $t>18$ .
  Besides  we have to take into account (\ref  {ExpressioneFratta}) which lowers further the contributions at  right. 
  Instead  spurious contribution at left are  multiplied by $\approx 2$ , but ,distance from $Boundary(L_1,L_2)    $, in worst case ($p=0$),  is, for $ t >20$: $ |a| = \sqrt{2\pi t} \approx \sqrt{t} \times  \sigma >  4.472 \sigma$.
  With ( \ref {ErfStuff}) we can appreciate the corresponding $\frac {\Delta R}{R}$  in  (\ref    {UpperBounsSum}) and  compare with data in
  fig. \ref {MaggiorazioniL023} .


\subsection {Upper bound of Main integral on $L_2$ path}

We follow 
  \cite[p.~142]{Edwards:1974cz}  `` Estimation of the integral away from saddle point'', second  case  j=2,with the only difference of introducing $\epsilon \ne 0$.
The boundary of   $L_2$ segment are:

\be \label{BoundaryL1L2}
Boundary(L_1,L_2) = a- \frac{|a|}{r} e^{i\frac{\pi}{4}} \quad ; \quad   1<r \le 2
\ee
$$
     Boundary(L_2,L_3) = ( Re[  a- \frac{|a|}{r} e^{i\frac{\pi}{4}}]  \quad  ,  \quad    -\pi(2N+1) ) \quad ; \quad   1<r \le 2
$$
 \noindent  Where  $a$ is  defined in (\ref {a}).
 
  \noindent    So , in the integration on $ L_2$,  always  holds :$$
\Re[x]=- \frac{|a|}{\sqrt{2}r}=  - \frac{\sqrt{\pi t}}{r}= -b
$$

 \noindent   The integrand denominator of  (\ref{MainInt})  is :
\be \label {LowerBoundDen}
|e^x-1 | \ge | e^{-b}-1 | > 0.5  \quad for \quad  t>10  \quad  1<r \le 2
\ee
As $-N (-b)= \left(-\sqrt{\frac{t}{2 \pi}}+p\right) (-b) < \sqrt{\frac{t}{2 \pi}} b$, the integrand numerator of  \ref{MainInt}  is at most:
$$
| (-x)^{\epsilon-\frac{1}{2}+it}e^{-Nx}| \le  max| x|^{\epsilon-\frac{1}{2}} max \quad e^{-t \ \Im[ln(-x)] } e^{\sqrt{\frac{t}{2 \pi}} b}
$$
Because  $\epsilon-1/2  \le 0$ the $ max| x|^{\epsilon-\frac{1}{2}} $ is positioned  where $L_2$ crosses  real axis, and,  it is: 
$$
 max| x|^{\epsilon-\frac{1}{2}} = b^{\epsilon-\frac{1}{2}}
$$
the  max value of  $ e^{-t  \  \Im[ln(-x)] }$ happens to be where the exponent is max, i.e. where $ \Im[ \ \ln(-x)]$  is min,  namely in point  $Boundary(L_1,L_2)$ where -x phase has a minimum.
We have ( $a$ is  defined in  \ref {a} ):
$$
 \Im[ \ \ln(-x)]  = \Im[ \ \ln(-\left( a -\frac{|a|}{r} e^{i\frac{\pi}{2} } \right) )]  = \Im \ln \left[  -i +\frac{1}{r \sqrt{2}}  +i\frac{1}{r \sqrt{2}}  \right] = \Im \ln \left[ (-i)\left(  1 -\frac{1}{r \sqrt{2}}  +i\frac{1}{r \sqrt{2}}  \right) \right] =
$$
$$
-\frac{\pi}{2}  + \arctan\left(    \frac{\frac{1}{r \sqrt{2}} }{1-\frac{1}{r \sqrt{2}} }         \right)= 
-\frac{\pi}{2}  + \arctan\left(    \frac{1}{\sqrt{2}r-1}        \right)
$$
 \noindent     So :
$$
max_{ on \  L_2} \quad \left\{   e^{-t  \ \Im[ \ ln(-x)] }  \right \} = e^{t\frac{\pi}{2}}e^{ -t \arctan \left(    \frac{1}{\sqrt{2}r-1}        \right)  }
$$
while
$$
e^{-N(-b)} < e^{\frac{t}{\sqrt{2}r}}
$$
Then the numerator  integrand module  in  (\ref{MainInt})  is at most:
$$
 b^{\epsilon-\frac{1}{2}}  e^{t\frac{\pi}{2}}e^{ -t \left ( 
 \arctan \left(    \frac{1}   {   \sqrt{2}r-1   } \right)   -         \frac{1}{\sqrt{2}r    }   
    \right)  }=
     b^{\epsilon-\frac{1}{2}}  e^{t\frac{\pi}{2}}e^{ -t \left ( 
 K_2(r)
    \right)  }
$$
Where:

\be \label {K2}
K_2(r) :=- \left(\arctan \left(    \frac{1}   {   \sqrt{2}r-1   } \right)   -         \frac{1}{\sqrt{2}r    } \right)  
\ee
So reminding that (\ref  {LowerBoundDen} ):
$$
\left(\frac{1}{e^x -1}\right)_{L_2} <2
$$
The length of segment  $L_2$  is about  $2|a|= 2 \sqrt{2 \pi t}= 2 r \sqrt{2}  \frac{\sqrt{\pi t}}{r} = 2r \sqrt{2} b$,
then, after clearing the common factor $F(t)$  and for  $ t >10$ ( see also ( \ref {MainInt} ) ) we have:
$$
\left|  e^{Re_2^-+Re_3^- )}     \frac{ e^{-\pi t/2}            }{   (  2 \pi  )^{\frac{1}{2} +\epsilon } (e^{-i\frac{\pi}{2} \epsilon}-i e^{i\frac{\pi}{2} \epsilon} e^{-t \pi}   )     } 
\left ( \sqrt{  \frac{2 \pi }{t}  } \right)^{\epsilon} 
\int_{L_2} \frac{(-x)^{s-1}e^{-Nx}dx}{e^x-1}   \right| < 
$$
$$
< \frac{ e^{-\pi t/2}            }{   (  2 \pi  )^{\frac{1}{2} +\epsilon }    } 
\left ( \sqrt{  \frac{2 \pi }{t}  } \right)^{\epsilon}   
 b^{\epsilon-\frac{1}{2}}  e^{t\frac{\pi}{2}}e^{ K_2(r) t  }    2r \sqrt{2} b=
$$
$$
= \frac{   b^{\epsilon-\frac{1}{2}}          }{   (  2 \pi  )^{\frac{1}{2} +\epsilon }    } 
\left ( \sqrt{  \frac{2 \pi }{t}  } \right)^{\epsilon}   
       e^{ K_2(r) t  }     2r \sqrt{2} b=  
    \frac{   \left(  \frac{\sqrt{\pi t}}{r}    \right)^{\epsilon-\frac{1}{2}}          }{   (  2 \pi  )^{\frac{1}{2} +\epsilon }    } 
\left ( \sqrt{  \frac{2 \pi }{t}  } \right)^{\epsilon}   
       e^{ -t \left ( 
 k
    \right)  }     2r \sqrt{2} \left(   \frac{\sqrt{\pi t}}{r}\right) = 
    $$
    
    $$
  . . . =     \frac{   \left(  \frac{\sqrt{\pi t}}{r}    \right)^{\epsilon}          }{   (  2 \pi  )^{\frac{1}{2} +\epsilon }    } 
\left ( \sqrt{  \frac{2 \pi }{t}  } \right)^{\epsilon}   
       e^{ K_2(r) t }     2r \sqrt{2} \left(   \frac{\sqrt{\pi t}}{r}\right)^{\frac{1}{2}} =
  \frac{   \left(       \frac{   1    }{r}    \right)^{\epsilon}          }{   (  2 \pi    )^{\frac{1}{2} +\epsilon }    } 
\left ( \sqrt{  \frac{2 \pi  \pi t }{t}  } \right)^{\epsilon}   
       e^{ -t \left ( 
 k
    \right)  }     2r \sqrt{2} \left(   \frac{\sqrt{\pi t}}{r}\right)^{\frac{1}{2}} =    . . .
$$

$$
= \frac{   \left(       \frac{   1    }{r}    \right)^{\epsilon}          }{   (  2 \pi    )^{\frac{1}{2} +\epsilon }    } 
\left ( \sqrt{   2} \pi  \right)^{\epsilon}   
       e^{K_2(r) t  }     2r \sqrt{2} \left(   \frac{\sqrt{\pi t}}{r}\right)^{\frac{1}{2}} =   ...
  \frac{   \left(       \frac{   1    }{r}    \right)^{\epsilon}          }{   (  2     )^{\frac{1}{2} +\epsilon }   \pi^{\frac{1}{2} } } 
\left ( \sqrt{   2}  \right)^{\epsilon}   
       e^{ -t \left ( 
 k
    \right)  }     2r \sqrt{2} \left(   \frac{\sqrt{\pi t}}{r}\right)^{\frac{1}{2}}   
$$

$$
 = \frac{   \left(       \frac{   1    }{r}    \right)^{\epsilon}          }{   (  2     )^{\frac{1}{2} +\epsilon }   \pi^{\frac{1}{2} } } 
\left ( \sqrt{   2}  \right)^{\epsilon}   
       e^{K_2(r) t }     2 \sqrt{r} \sqrt{2} \left(   \sqrt{\pi t}    \right)^{\frac{1}{2}}   =
  \frac{   \left(       \frac{   1    }{r\sqrt{   2}}    \right)^{\epsilon}          }{   (  2     )^{\frac{1}{2} }   \pi^{\frac{1}{2} } }  
       e^{ -t \left ( 
 k
    \right)  }     2 \sqrt{r} \sqrt{2} \left(   \sqrt{\pi t}    \right)^{\frac{1}{2}}   = . . .
$$

\be \label {MaggL2}
 . . . =
  \frac{   \left(       \frac{   1    }{r\sqrt{   2}}    \right)^{\epsilon}          }{     \pi^{\frac{1}{2} } }  
       e^{ K_2(r) t  }     2 \sqrt{r}  \left(   \sqrt{\pi t}    \right)^{\frac{1}{2}}   = 
     \frac{   \left(       \frac{   1    }{r\sqrt{   2}}    \right)^{\epsilon}          }{     \pi^{\frac{1}{4} } }  
       e^{ K_2(r) t }     2 \sqrt{r}  \left(   t \right)^{\frac{1}{4}} =        
\ee
$$   \left(  \frac{   1    }{r\sqrt{   2}}    \right)^{\epsilon}  e^{K_2(r) t }     2 \sqrt{r}  \left(   t/\pi \right)^{\frac{1}{4}} =  UpperBound_{L_2}(t)$$
Note that for $\epsilon >0 $  the upper bound (\ref {MaggL2} ) decreases.

\subsection {Upperbound of Main integral on $L_3$ path}

We follow  \cite[p.~144]{Edwards:1974cz}   `` Estimation of the integral away from saddle point'', third case  j=3 ,with the only difference of introducing $\epsilon > 0$.

$$
Boundary(L_2,L_3)= \left(-b \quad,\quad -\pi i ( 2N+1) \right)=\left(-\frac{\sqrt{\pi t}}{r}\quad,\quad -\pi i ( 2N+1) \right)
$$
Denominator  is :
$$
e^x-1 = e^{ \Re(x) - \pi i ( 2N+1)}-1 = -e^{ \Re(x) }-1 \rightarrow | -e^{ \Re(x) }-1| >1
$$
so module of denominator  is always $>1$.
We can write integral numerator as :
\be \label {IntegrNumerator}
(-x)^{\epsilon-1/2} (-x)^{it} e^{-Nx}
\ee
Because $\epsilon-0.5 \le 0$  the  max  of  $|x|^{\epsilon-1/2}$, on $L_3$,  is:  $|\pi ( 2N+1)|^{\epsilon-1/2}$, attained when $\Re(x)=0$.

$$
|(-x)^{it} |=|e^{(\Re[ \ \ln(-x)] + i \Im[ \ \ln(-x)]) it}| =e^{-t  \ \Im[ln(-x)]}
$$
What about max  reached by $-t  \ \Im[ \ \ln(-x)]$  on $L_3$?

\noindent  Namely which is the smallest phase  of -x on  $L_3$?.

\noindent  Note that $\angle x = \angle [-(-x )] $  is negative and then increases , so $\angle [-x ]$ is positive in $Boundary(L_2,L_3)$  and then decreases  (i.e.  $\Re(x) \rightarrow + \infty$ namely $\Re(-x) \rightarrow - \infty$ ).

 \noindent Answer: the phase assumed in $Boundary(L_2,L_3)$:
$$
\arctan\left( \frac{\pi(2N+1)}{b}  \right ) \approx \arctan\left( \frac{\pi \left(2\sqrt{\frac{t}{2 \pi}}+1 \right)}{\frac{\sqrt{t 2 \pi}}{r}}  \right ) >\frac{\pi}{4}
$$
because  $| b| =\frac {\sqrt{t \pi}}{r} $ and:
$$
\frac{\pi \left(2\sqrt{\frac{t}{2 \pi}}+1 \right)}{\frac{\sqrt{t 2 \pi}}{r}} =\frac{r}{\sqrt{t \pi}  } \left(  \sqrt {2 t \pi} +\pi \right) = r \left (   \sqrt(2) +\sqrt{\frac{\pi}{t}} \right) >1  \quad ; because \quad  1 <r \le2
$$
So 

   \be \label {BoundOnNumer}
   |(-x)^{it} |< e^{-t \frac{\pi}{4}}    
    \ee
So the module of \ref{MainInt} on  $L_3 \quad  (t>10)$, putting $u= \Re(x)$  and  taking off  $F(t)$  factor, 
for \ref  {BoundOnNumer}  is less then:
$$
\left[    \left( \frac{1}{2 \pi t}  \right)^{\epsilon/2}        \frac{    e^{   -t \pi   /   2}    } {   \sqrt{2 \pi}      }     \right]
e^{-t \pi/4} [ (2N+1)\pi]^{\epsilon-1/2} \int_{-b}^\infty e^{-Nu} du  
$$
\be \label {MaggL3}
 =\left(   \frac{\sqrt{2 \pi t}+\pi}{\sqrt{2 \pi t}} \right)^{\epsilon} \frac{ e^ {   -3t \pi/4  }  }   {  \sqrt{2}\pi \sqrt {2N+1}   }    \frac{    e^{     t   /    (   \sqrt{2}  r)     }  }   { N  } < \left(   \frac{\sqrt{2 \pi t}+\pi}{\sqrt{2 \pi t}} \right)^{\epsilon}  \frac {e^{-t}} { \sqrt{\pi t}} 
 =  UpperBound_{L_3}(t)
\ee
because  
$$\forall t>10 \quad    (e^{-i\frac{\pi}{2} \epsilon}-i e^{i\frac{\pi}{2} \epsilon} e^{-t \pi}   )\approx e^{-i\frac{\pi}{2} \epsilon}$$
 And considering that
$$\left| \frac {e^{i\frac{\pi}{4}\epsilon} e^{-i\theta_1(t)}}{(2 \pi  )^{ +it} e^{-i\frac{\pi}{4}}    e^{-i\frac{\pi}{2} \epsilon}   } \right|  =1$$
then we have:
$$
\left|
\frac{ e^{i\frac{\pi}{4}\epsilon} e^{-\pi t/2} e^{-i\theta_1(t)} }{   (  2 \pi  )^{\frac{1}{2} +\epsilon +it} (e^{-i\frac{\pi}{2} \epsilon}-i e^{i\frac{\pi}{2} \epsilon} e^{-t \pi}   ) e^{-i\frac{\pi}{4}}   } 
\left ( \sqrt{  \frac{2 \pi }{t}  } \right)^{\epsilon} \right|=
\frac{ e^{-\pi t/2}  }
{   (  2 \pi  )^{\frac{1}{2} +\epsilon } 
  } 
\left ( \sqrt{  \frac{2 \pi }{t}  } \right)^{\epsilon} 
 =
\left[    \left( \frac{1}{2 \pi t}  \right)^{\epsilon/2}        \frac{    e^{   -t \pi   /   2}    } {   \sqrt{2 \pi}      }     \right]
$$
While, for \ref  {BoundOnNumer} and \ref  {IntegrNumerator}:

$$
\left|   \int_{L_3}  \frac{   (-x)^{-\frac{1}{2} +\epsilon +it }  e^{-Nx} dx }{ e^x-1}   \right| <e^{-t \pi/4} [ (2N+1)\pi]^{\epsilon-1/2} \int_{-b}^\infty e^{-Nu} du
$$
besides:

\be \label {NPerMenB}
N [-(-b)] \approx  \sqrt{ \frac {t}{2 \pi}}  \frac {\sqrt{t \pi}}{r} =  \frac {   t   }    {   \sqrt{2}  r}   
\ee
and:

$$
  -3t \pi/4 +  t   /    (   \sqrt{2}  r)  =t \left(- \frac {3 \pi}{4}  +\frac {1}{\sqrt{2}r}\right) < < -t
$$
Note that  upper bound on $L_3$  \ref {MaggL3} $ e^{-t}$ increases with $\epsilon$ but  it is irrelevant in front of  \ref {MaggL0}, and \ref  {MaggL2}, see fig   \ref {MaggiorazioniL023}.

\subsection {Conclusions on Upperbounds }                                                                                                                                                                                                                                                                      \label {MaggiorazioniL023}

Upper bounds can be summed up   and compared with minimum value of $R_0(t)  $ ( i.e. when $p=0.5$, see (\ref {CZero}) and (\ref {R0DiT})   ) . 

\noindent Summarizing results we have:
      
\noindent (\ref{MaggL0}) (i.e. (\ref{MainInt}) along  $L_0$ upperbound )  is decreasing if $\epsilon$ grows.
  $K_0(r)$ is defined in  \ref {K0},   $K_0(r) \le K_0(2) \approx -0.1 \  \ 1 <r \le 2$. 
 
 \be \label {UBL0}
\frac{    \int_{L_0} \frac{(-x)^{s-1}e^{-Nx}dx}{e^x-1}   }{ R_0(t)_{p=0.5} }
< \frac {
 \frac { 2 \sqrt{2}    }{    \left[\sqrt{2 \pi } \right]^{ (1+\epsilon )/2 }     \sqrt{ \sqrt { 2 \pi t   } } }
 }
 { R_0(t)_{p=0.5}}  \ e^{K_0(r) t}= UpperBound_{L_0}(t)
 \ee
 \noindent    \ref{MaggL2} (i.e. \ref{MainInt} along  $L_2$ upperbound) is decreasing if $\epsilon$ grows.
  $K_2(r)$ is defined in  \ref {K2},  $K_2(r) \le K_2(2) \approx -0.2 \  \ 1 <r \le 2$. 
  \be \label{UBL2}
\frac{    \int_{L_2} \frac{(-x)^{s-1}e^{-Nx}dx}{e^x-1}   }{ R_0(t)_{p=0.5} }
< \frac {
       \left(  \frac{   1    }{r\sqrt{   2}}    \right)^{\epsilon}          
           2 \sqrt{r}  \left(   t/\pi \right)^{\frac{1}{4}} 
 }
 { R_0(t)_{p=0.5}}  \ e^{K_2(r) t }  =  UpperBound_{L_2}(t)
 \ee
 
\begin{figure}[!htbp]
\begin{center}
\includegraphics[width=1.0\textwidth]{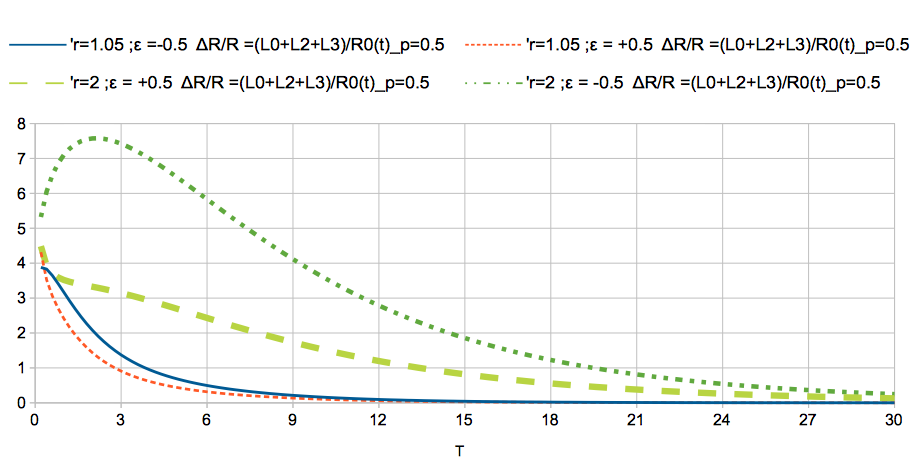} 
\caption{\small { 
 Comparison 
 of  overall Upperbounds ratio $\frac{\Delta R}{R}$   (\ref  {UpperBounsSum})  with parameters $r$ (\ref {BoundaryL0L1}) and $\epsilon$. With   $(r,\epsilon) = (1.05, \pm 0.5)$   we have the two  lower curves,and,  $\epsilon= +0.5$ is the lower one between the two.. For upper curves we have  $(r,\epsilon) = (2, \pm 0.5)$, and, the lower curve is always with  $\epsilon= +0.5$.Of course the case $\epsilon=0$ is  midway between  found bounds for  $\epsilon=\pm 1/2$.  Due to the conjugate symmetry of $\xi(1/2+\epsilon+it)$ with respect to $\epsilon$, it is the lower UpperBound   that applies. So errors in(\ref {espressione1}) ,and similar formulas, is maximum for $\epsilon=0$. {\bf The aim of error upperbounds in extended  $Z(t,\epsilon)$ less  or equal to the same in $Z(t)$  is so  met}
Here , differently as  in  \cite[p.~141]{Edwards:1974cz}  ,  (\ref {BoundaryL0L1}  )and   (\ref {BoundaryL1L2} ) , in $     a \pm \frac{|a|}{r} e^{i\frac{\pi}{4}} \quad  $  the parameter $r$ has been left free to change  in the interval $ \quad   1<r \le 2 $  .  This choice gives a more realistic evaluation  of upperbounds error of ( \ref {espressione1} )  ,and similar formulas, at low t values. In   \cite[p.~144]{Edwards:1974cz}     it is suggested  the value    $t \ge 100$ to get a practically  error free  expression . Considering the figure, it  appears pessimistic, and ,   only with $r$ close to 2 ( as chosen in  \cite  {Edwards:1974cz})  is justified .  The only restriction on $r$ is linked to convergence of  (\ref {LogDevelopment} )       .  The convergence speed is not an issue.  By the way this explains  the good match noted in     \cite[p.~155]{Edwards:1974cz}        between  $ Z(t>18,\epsilon=0)$ and Haselgrove table      \cite[p.~122]{Edwards:1974cz}    . 
} }
\label {MaggiorazioniL023}
\end{center}
\end{figure}
 \noindent    \ref{MaggL3} (i.e. \ref{MainInt} along  $L_3$ upperbound  ) increases with $\epsilon$ but  the factor $e^{-t}$ cancel the contribution .
  \be \label{UBL3}
\frac{    \int_{L_3} \frac{(-x)^{s-1}e^{-Nx}dx}{e^x-1}   }{ R_0(t)_{p=0.5} }
< \frac {
       \left(   \frac{\sqrt{2 \pi t}+\pi}{\sqrt{2 \pi t}} \right)^{\epsilon} }
 { R_0(t)_{p=0.5}}   \  \frac {e^{-t}} { \sqrt{\pi t}  }= UpperBound_{L_3}(t)
 \ee
Of course the upper bound of the sum is less then the sum of the upper bounds ,
  so we have an overall  upperbound ratio of:
  
  \be  \label {UpperBounsSum}
  \left|
\frac{    \int_{L_0,L_2,L_3} \frac{(-x)^{s-1}e^{-Nx}dx}{e^x-1}   }{ R_0(t)_{p=0.5} }
\right|
= \frac {\Delta R}{R}
<UpperBound_{L_0}(t)+UpperBound_{L_2}(t)+UpperBound_{L_3}(t)
\ee
Of course the case $\epsilon=0$ is between $\epsilon=\pm 1/2$ .
Here, in (\ref {BoundaryL0L1}) and in (\ref {BoundaryL1L2}), differently that  in     \cite[p.~141]{Edwards:1974cz}   
the parameter $r$ has been left free to change  in the interval $ \quad   1<r \le 2 $  .  This choice gives a more realistic evaluation  of upper-bounds error of (\ref {espressione1})  and similar formulas at low t values. In      \cite[p.~144]{Edwards:1974cz}          it is suggested  the value    $t \ge 100$ to get a practically  error free  expression.\\
In light of  above  consideration, this  appears pessimistic because   only with $r$ close to 2 is justified.   The only restriction on $r$ is linked to convergence of  (\ref {LogDevelopment}).  The convergence speed is not an issue.  

\noindent    By the way, this explains  the good match noted in   \cite[p.~155]{Edwards:1974cz}     between  $ Z(t>18,\epsilon=0)$ and Haselgrove table      \cite[p.~122]{Edwards:1974cz}. 
The spurious contribution  (\ref {GammaRef}) is  an insignificant fraction of  $R_0(t)$  at least for $t>20$.
 
\noindent Besides  the bound with   $\epsilon >0$   ( and so also for $\epsilon <0$)  is lower  than with  $\epsilon =0$.

\noindent We  can affirm that   (\ref  {ZSinhECosh}) 
 can be used profitably from t =20 with $\epsilon$ in critical strip and beyond.


\section*{Acknowledgments}

I thank professor Richard B. Paris for his kind support and for useful discussions.
I also thank Paolo Lodone for useful discussions and for contributing in some points of this work.

\section*{Bibliography}

\vspace{0.3cm}



\end{document}